
\magnification = \magstep1
\vsize = 9truein

\def\Z{\ifx\MYUN\Bbb {\font\bb=msbm10 \hbox{\bb Z}}
\else {\Bbb Z} \fi}
\def\C{\ifx\MYUN\Bbb {\font\bb=msbm10 \hbox{\bb C}}
\else {\Bbb C} \fi}
\def\R{\ifx\MYUN\Bbb {\font\bb=msbm10 \hbox{\bb R}}
\else {\Bbb R} \fi}
\def\Q{\ifx\MYUN\Bbb {\font\bb=msbm10 \hbox{\bb Q}}
\else {\Bbb Q} \fi}
\def\qed{{~~~\vrule height .75em width .4em depth .3em}}
\def\irr#1{{\rm Irr}(#1)}

\def\deg#1{{\rm deg}(#1)}
\def\nor{\triangleleft\,}
\def\gen#1{\langle#1\rangle}
\def\zent#1{{\bf Z}(#1)}
\def\ker#1{{\rm ker}(#1)}
\def\ritem#1{\item{{\rm #1}}}
\def\iitem#1{\goodbreak\par\noindent{\bf #1}}
\def\inv{^{-1}}
\def\aut#1{{\rm Aut}(#1)}
\def\plusdot{\buildrel\textstyle.\over+}

\def\sbs{\subseteq}
\def\gcd#1#2{{\rm gcd}(#1,#2)}
\def\dimm#1#2{{\rm dim}_{#1}(#2)}
\def\tens{\otimes}

\def\th{\,{\rm th}}
\def\mdot{\hbox{$\cdot$}}
\def\ref#1{{\bf[#1]}}

\def\bur{1}
\def\camtt{2}
\def\cam{3}
\def\donsta{4}
\def\ie{5}
\def\hup{6}
\def\isa{7}
\def\isc{8}
\def\rud{9}
\def\stelen{10}
\def\tao{11}
\def\ter{12}
\def\wie{13}

\def\supp#1{{\rm supp}(#1)}
\def\XX{{\cal X}}
\def\BB{{\cal B}}
\def\RR{{\cal R}}
\def\TT{{\cal T}}

{\nopagenumbers
\font\big=cmbx12

\vglue .2truein
\big
\centerline{INEQUALITIES FOR}
\smallskip
\centerline{FINITE GROUP}
\smallskip
\centerline{PERMUTATION MODULES}
\bigskip\rm
\centerline{by}
\bigskip
\bf
\centerline{Daniel Goldstein}
\centerline{Center for Communications Research}
\centerline{4320 Westerra Ct.}
\centerline{San Diego CA~~92121}
\smallskip
\centerline{E-Mail: dgoldste@ccrwest.org}
\bigskip

\centerline{Robert M. Guralnick}
\centerline{Department of Mathematics}
\centerline{University of Southern California}
\centerline{1042 W. 36th Place}
\centerline{Los Angeles CA~~90089}
\smallskip
\centerline{E-Mail: guralnic@math.usc.edu}
\bigskip

\centerline{I.~M.~Isaacs}
\centerline{Mathematics Department}
\centerline{University of Wisconsin}
\centerline{480 Lincoln Drive}
\centerline{Madison WI~~53706}
\smallskip
\centerline{E-Mail: isaacs@math.wisc.edu}
\bigskip
\rm
\vfil

{\narrower\noindent
{\bf ABSTRACT:}~~If $f$ is a nonzero complex-valued function
defined on a finite abelian group $A$ and $\hat f$ is its Fourier
transform, then  $|\supp f||\supp{\hat f}| \ge |A|$, where
$\supp f$ and $\supp{\hat f}$ are the supports of $f$ and $\hat f$. In
this paper we generalize this known result in several directions. In
particular, we prove an analogous inequality where the abelian
group $A$ is replaced by a transitive right $G$-set, where $G$ is an
arbitrary finite group. We obtain stronger inequalities when the
$G$-set is primitive and we determine the primitive groups for which
equality holds. We also explore connections between inequalities of
this type and a result of Chebotar\"ev on complex roots of unity, and
we thereby obtain a new proof of Chebotar\"ev's theorem.\par}

\vfil

\font\sf=cmr7
\sf\noindent
The research of the second author was partially supported by
Grant DMS 0140578 of the U.~S.~NSF.\par
\noindent
The research of the third author was partially supported by the
U.~S.~NSA. 
\eject}

\iitem{1.~Introduction.}

The starting point for this paper is an inequality for complex-valued
functions defined on finite abelian groups. This result is generally
attributed to D.~L.~Donoho and P.~B.~Stark, although in their
paper \ref\donsta, they prove the result only for cyclic groups.
(The more general statement, with a proof, can be found as
Theorem~14.1 in the book by A.~Terras \ref\ter.)
\medskip

\iitem{THEOREM  A (Donoho and Stark).}~~{\sl Let $A$ be a finite
abelian group and suppose that $f$ is an arbitrary nonzero
complex-valued function on $A$. Writing $\hat f$ to denote the
Fourier transform of $f$, we have
$$
|\supp f||\supp{\hat f}| \ge |A| \,.
$$}
\medskip

Here, $\supp f$ is the {\bf support} of $f$, which is the set of
elements of $A$ on which the function $f$ takes nonzero values.
Similarly, $\supp{\hat f}$ is the set of linear characters of $A$ on
which $\hat f$ takes nonzero values. In other words, a linear
character $\lambda$ of $A$ lies in $\supp{\hat f}$ precisely when
$\lambda$ occurs with nonzero coefficient when the function $f$ is
written as a linear combination of the set $\hat A$ of linear
characters of $A$. 

We can think about Theorem~A in the following much more general
context. Let $S$ be a finite right $G$-set, where $G$ is an arbitrary
group, and suppose that the action of $G$ on $S$ is transitive. (It
will be no loss to assume that the action of $G$ on $S$ is faithful,
and so we can assume that $G$ is finite.) Let $F$ be an arbitrary
field and write $F[S]$ to denote the $F$-space of $F$-linear
combinations of members of $S$. (Alternatively, we could view this
space as the set of $F$-valued functions on $S$, but we prefer the
former point of view, in which $S$ is a subset of $F[S]$.) If we
extend the given action of $G$ on $S$ linearly to all of $F[S]$, then
$F[S]$ becomes a right $G$-module over $F$: the permutation
module. Now if $v \in F[S]$, we write $\supp v$ to denote the set
of points of $S$ that occur with nonzero coefficients in $v$. (This,
of course, is exactly the support of the function $f$ corresponding
to the vector $v$, where $v$ and $f$ are related by the equation
$v = \sum_{s \in S} f(s)s$.)

Returning now to the situation of Theorem~A, take $F = \C$, the
complex numbers, and $G = S = A$, where the transitive action of
$G$ on $S$ is the regular action, defined by multiplication in $A$.
In the theorem, we are given a function $f$ defined on $A$, and
we let $v \in F[S] = \C[A]$ be the corresponding vector, so that
$\supp f = \supp v$. But how can we interpret $\supp{\hat f}$ in
terms of the vector $v$?

View $F[S] = \C[A]$ as the complex group algebra of $A$, and
write $e_\lambda \in \C[A]$ to denote the idempotent
corresponding to the linear character $\lambda$ of $A$. (Recall
that $e_\lambda = (1/|A|) \sum_{a \in A} \overline{\lambda(a)} a$.)
Because $A$ is abelian, the idempotents $e_\lambda$ form a basis
for $\C[A]$, and thus there exist coefficients $a_\lambda \in \C$
such that $v = \sum_\lambda a_\lambda e_\lambda$. If $\lambda$
and $\mu$ are distinct linear characters of $A$ and we view them
(by linear extension) as being defined on the entire group algebra
$\C[A]$, then $\lambda(e_\lambda) = 1$ and $\mu(e_\lambda) = 0$.
It follows that
$$
a_\lambda = \lambda(v) = \sum_{g \in A} \lambda(g)f(g)
= |A| \hat f (\overline\lambda) \,,
$$
where the third equality follows from the definition of $\hat f$. We
see, therefore, that $|\supp{\hat f}|$ is equal to the number of
coefficients $a_\lambda$ that are nonzero. 

It is well known that the ideals of $\C[A]$ are exactly the subspaces
of $\C[A]$ spanned by the various subsets of the set
$\{e_\lambda \mid \lambda \in \hat A\}$. (This, of course, is because
the group $A$ is abelian. The corresponding general statement for
an arbitrary finite group $G$ is that the ideals of the center
$\zent{\C[G]}$ of the group algebra are the subspaces spanned by
the central idempotents $e_\chi$, where $\chi$ runs over $\irr G$,
the set of irreducible characters of $G$.)

It follows that $|\supp{\hat f}|$ is the dimension of the smallest
ideal of $\C[A]$ that contains $v$. Since in this case, the ideals
of $\C[A]$ are precisely the $G$-submodules of $F[S]$, we
conclude that $|\supp{\hat f}|$ is equal to the dimension of the
$G$-submodule $\gen{v}$ of $F[S]$ generated by $v$. We see,
therefore, that the following result includes Theorem~A.
\medskip

\iitem{THEOREM B.}~~{\sl Let $S$ be a finite transitive right $G$-set
and let $0 \ne v \in F[S]$, where $F$ is an arbitrary field. Then
$$
|\supp v|\dim(\gen v) \ge |S| \,,
$$
where $\gen v$ is the $G$-submodule of $F[S]$ generated by $v$.}
\medskip

Observe that this result generalizes Theorem~A in three distinct
ways. The group $G$ need not be abelian and the field $F$ need
not be the complex numbers. Also, the function $f$ need not be
defined on the group $G$ itself; it can be defined instead on an
arbitrary transitive right $G$-set $S$.

Our proof of Theorem~B is quite easy, and in fact, it is shorter than
the proof of Theorem~A given in \ref\ter.
\medskip

\iitem{Proof of Theorem~B.}~~The submodule $\gen v \sbs F[S]$ is
the linear span of the $G$-translates $vg$ of $v$ as $g$ runs over
$G$, and so we can choose a basis $\BB$ for $\gen v$ consisting of
such translates. Since $v$ is nonzero and $G$ acts transitively on
$S$, we see that every point of $S$ is in $\supp{vg}$ for some
element $g \in G$. But $vg$ is a linear combination of members of
the basis $\BB$, and thus $S = \bigcup_{b \in \BB} \supp b$. It
follows that $\sum_{b \in \BB} |\supp b| \ge |S|$. If $b \in \BB$,
however, then $b$ is a translate of $v$, and we see that
$|\supp b| = |\supp v|$. We conclude that $|\BB||\supp v| \ge |S|$,
and since $|\BB| = \dim(\gen v)$, the proof is complete.\qed
\medskip

In order to facilitate discussion of Theorem~B and related results,
we establish the following notation. Given $v \in F[S]$, write
$t = t(v) = |\supp v|$ and $d = d(v) = \dimm F{\gen v}$, and set
$n = |S|$. The assertion of Theorem~B, therefore, is that if $t > 0$,
then $td \ge n$.

It is not hard to describe exactly when equality holds in Theorem~B,
and we discuss this next. It is clear from our proof that if equality
holds, then the supports of distinct members of the basis $\BB$ must
be disjoint. Since $\BB$ can be chosen to contain any two linearly
independent translates of $v$, it follows that if $\Delta = \supp v$
and $g \in G$, then either $\Delta g \cap \Delta = \emptyset$, or else
$vg = av$ for some nonzero scalar $a \in F$, and in particular,
$\Delta g = \Delta$. 

Recall that a nonempty subset $\Delta$ of the $G$-set $S$ called
a {\bf block} if for each element $g \in G$, either
$\Delta g = \Delta$ or $\Delta g \cap \Delta = \emptyset$. Of course,
$\Delta = S$ is a block, and so is $\Delta = \{s\}$ for each point
$s \in S$. These are the {\bf trivial} blocks, and we recall that the
transitive $G$-set $S$ is said to be {\bf primitive} if every block is
trivial. Observe that if $\Delta$ is a block, then so is every
$G$-translate, and thus the distinct translates of $\Delta$ partition
$S$. The number of such translates, therefore, is $|S|/|\Delta|$.

If equality holds in Theorem~B, we know that $\Delta = \supp v$
must be a block. Conversely, given any block $\Delta \sbs S$,
consider the vector $v = \sum_{x \in \Delta} x$. Then
$\supp v = \Delta$, and in fact, the supports of the distinct
$G$-translates of $v$ are exactly the $|S|/|\Delta| = n/t$ distinct
translates of $\Delta$. The translates of $v$ are thus linearly
independent and we have $d = d(v) = n/t$, and so equality holds in
Theorem~B. 

It is possible to have $td = n$ even if $v$ is not a scalar multiple of
the sum of the points in its support. We digress briefly to discuss the
most general possible case where this equality holds. Assuming
equality, write $\Delta = \supp v$ and let $H = G_\Delta$ be the
(setwise) stabilizer of the block $\Delta$. If $h \in H$, then the
translate $vh$ also has support $\Delta$, and thus, as we have seen,
$vh$ must be a scalar multiple of $v$. If we write
$vh = \lambda(h) v$, where $\lambda(h) \in F$, it is easy to see that
$\lambda$ is a homomorphism from $H$ into $F^\times$. Also, if
$x \in \Delta$ and $K = G_x$ is the stabilizer of $x$, then $K \sbs H$
since $\Delta$ is a block. If $k \in K$, we must have $\lambda(k) = 1$
because the (nonzero) coefficients of $x$ in $v$ and in
$vk = \lambda(k) v$ are equal. 

In the situation of the previous paragraph, the vector $v$ is uniquely
determined (up to a scalar multiple) by the block $\Delta$ and the
homomorphism $\lambda$ from $H = G_\Delta$ into $F^\times$. To
see this, fix $x \in \Delta$ and note that $xH = \Delta$. If the
coefficient of $x$ in $v$ is $a$, then the coefficient of $xh$ in $v$
equals the coefficient of $x$ in $vh\inv = \lambda(h\inv)v$. Since
this coefficient is $a\lambda(h\inv)$, it follows that $v$ is
determined by $\Delta$, $\lambda$ and the scalar $a$, as claimed.

Conversely, let $\Delta$ be any block of the right $G$-set $S$. Write
$H = G_\Delta$ and let $\lambda: H \to F^\times$ be a homomorphism
such that $G_x \sbs \ker\lambda$ for $x \in \Delta$. We construct a
corresponding vector $v \in F[S]$ as follows. Fix $x \in \Delta$ and
recall that $xH = \Delta$. Write $v = \sum_{y \in \Delta} a_y y$, where
$a_{xh} = \lambda(h\inv)$, and note that this is well defined because
$G_x \in \ker\lambda$. It is easy to check that the supports of the
translates of $v$ are exactly the $|S|/|\Delta|$ translates of
$\Delta$, and that any two translates of $v$ with equal supports are
scalar multiples of one another. It follows that
$d = d(v) = |S|/|\Delta|$, and thus $td = n$, as wanted. 

If $S$ is a primitive $G$-set, then all blocks are trivial, and so if
$1 < t < n$, then equality cannot hold in Theorem~B and we have
$td > n$. In fact, an even stronger inequality holds in the primitive
case.
\medskip

\iitem{THEOREM C.}~~{\sl Let $S$ be a finite primitive right $G$-set
and suppose that $v \in F[S]$, where $F$ is an arbitrary field. If
$1 \le t < n$, we have $(t+1)d \ge 2n$, where as usual,
$t = |\supp v|$, $d = \dimm F{\gen v}$ and $n = |S|$.}
\medskip

We shall see that equality can hold in Theorem~C, but only under
highly restrictive conditions.
\medskip

\iitem{THEOREM D.}~~{\sl Suppose that $(t+1)d = 2n$ in
Theorem~C, where $1 < t < n-1$. Then $t = d$ and the action of $G$
on the set $\Omega$ of $G$-translates of $\supp v$ is doubly
transitive. Furthermore, every point of $S$ lies in exactly two
members of $\Omega$ and every two distinct members of $\Omega$
intersect in a single point. Finally, the field $F$ must have
characteristic $2$ and $v$ is a scalar multiple of the sum of the
points in its support.}
\medskip

Note that in the situation of Theorem~D, the group $G$ is a doubly
transitive permutation group in which the action on the collection
of two-point subsets is primitive. Conversely, we shall see that
given any such doubly transitive group and any field of
characteristic $2$, there is a corresponding example where equality
holds in Theorem~C. In Section~4, we describe exactly which
doubly transitive groups act primitively on the two-point subsets. 

If $n = |S|$ is a prime number, then the $G$-set $S$ is
automatically primitive. In this case, and provided that the field
$F$ has characteristic $0$, we prove an even better inequality than
that of Theorem~C. 
\medskip

\iitem{THEOREM E.}~~{\sl Let $0 \ne v \in F[S]$, where $S$ is a
transitive right $G$-set of prime cardinality $n$ and $F$ has
characteristic $0$. Then $t + d > n$, where $t = |\supp v|$ and
$d = \dimm F{\gen v}$.}
\medskip

We show in Section~6 that Theorem~E is essentially equivalent to a
certain theorem of N.~G.~Chebotar\"ev concerning complex roots
of unity of prime order. We give an independent proof of
Theorem~E, and this, in turn, yields a new proof of Chebotar\"ev's
result.

While this paper was in its final stages of preparation, the authors
learned that what is essentially our Theorem~E was independently
(and approximately simultaneously) discovered by T.~Tao. His
preprint \ref\tao\ presents a proof of Chebotar\"ev's result and
deduces Theorem~E from it. Also, Tao credits A.~Bir\'o with
independent discovery of the same inequality.

The inequality $t + d > n$ of Theorem~E is the best we can ever
hope to prove under any set of hypotheses. This fact is a
consequence of the following easy lemma.
\medskip

\iitem{LEMMA F.}~~{\sl Let $S$ be a finite transitive right $G$-set
and suppose that $M \sbs F[S]$ is an arbitrary nonzero
$G$-submodule, where $F$ is an arbitrary field. Then there exists a
nonzero vector $v \in M$ such that $t + d \le n + 1$, where $t$, $d$
and $n$ have their usual meanings.}
\medskip

\iitem{Proof.}~~Since $n + 1 - \dim(M) \le n = |S|$, we can choose a
subset $X \sbs S$ such that $|X| = n + 1 - \dim(M)$. Let
$W \sbs F[S]$ be the space of all vectors with support contained in
$X$, and note that $\dim(W) = |X|$. We have
$\dim(W) + \dim(M) = n + 1 > \dim(F[S])$, and it follows that
$W \cap M > 0$. Let $v$ be a nonzero vector in this intersection.
Then $\supp v \sbs X$, and so
$t = |\supp v| \le |X| = n + 1 - \dim(M)$. Also, $\gen v \sbs M$, and
thus $d = \dim(\gen v) \le \dim(M)$ and we have
$t + d \le (n + 1 - \dim(M)) + \dim(M) = n + 1$, as wanted.\qed
\medskip

We close this introduction by mentioning that there already is in
the literature a generalization to nonabelian groups of the theorem
of Donoho and Stark. (See Section 8 of \ref\donsta.) This result,
which we paraphrase somewhat, is due to P.~Diaconis and
M.~Shahshahani. It is much weaker than our Theorem~B, and it
follows as a corollary of our result.
\medskip

\iitem{COROLLARY G.}~~{\sl Let $f$ be a nonzero complex-valued
function on an arbitrary finite group $G$. For each character
$\chi \in \irr G$, choose a complex representation $\RR$ that
affords $\chi$, and construct the matrix
$M_\chi = \sum_{g \in G} f(g)\RR(g)$. Let $\XX$ be the
set of characters $\chi \in \irr G$ such that $M_\chi$ 
is not the zero matrix. Then
$$
|\supp f|\sum_{\chi \in \XX} \chi(1)^2 \ge |G| \,.
$$}
\medskip

Of course, the representation $\RR$ is determined by the character
$\chi$ only up to similarity, but this ambiguity does not affect
whether or not $M_\chi = 0$. Note that if $G$ is abelian, then the
set $\XX$ consists exactly of the complex conjugates of the linear
characters in the support of the Fourier transform $\hat f$. We
see, therefore, that Corollary~G reduces to Theorem~A in the
abelian case.
\medskip

\iitem{Proof of Corollary~G.}~~Let $v \in \C[G]$ be the vector
corresponding to the function $f$. Recall that $\C[G]$ is the direct
sum of its minimal ideals, and these ideals correspond to the
irreducible characters of $G$. If we extend the representation
$\RR$ affording $\chi$ to the whole group algebra by linearity, then
$\RR$ is exactly the projection map of $\C[G]$ onto the direct
summand corresponding to $\chi$. Also, $M_\chi = \RR(v)$, and so it
follows that $v$ lies in the sum of the minimal ideals for which
$M_\chi \ne 0$. This sum has dimension
$\sum_{\chi \in \XX} \chi(1)^2$, and it clearly contains the right
ideal of $\C[G]$ generated by $v$. Since $\supp f = \supp v$, we
see now that the result follows by Theorem~B applied to the action
of $G$ on itself by right multiplication.\qed
\bigskip

\iitem{2. Primitive actions.}

In this section we prove Theorem~C, which asserts the inequality
$(t+1)d \ge 2n$ in the case where $S$ is a primitive $G$-set and
$1 \le t < n$. Also in this section we study the case where
$(t+1)d = 2n$ and we prove Theorem~D.

The key to proving the inequality of Theorem~C is Rudio's lemma
\ref\rud, which also appears as 8.2 in Chapter~I of \ref\wie.
Actually, we need only the following weak form of Rudio's lemma,
which is valid even if $S$ is infinite. (Of course, for the purposes of
this paper, we need only the finite case.) 
\medskip

\iitem{(2.1) LEMMA.}~~{\sl Let $S$ be a primitive right $G$-set and
let $X$ be an arbitrary nonempty proper subset of $S$. If
$u,v \in S$ are distinct, then there exists some $G$-translate of $X$
that contains exactly one of $u$ and $v$.}
\medskip

The full statement of Rudio's lemma asserts that if $G$ is finite, we
can prescribe which one of $u$ or $v$ lies in a translate of $X$. As
we shall see in Section~7, however, that conclusion does not
necessarily hold for infinite groups.
\medskip

\iitem{Proof of Lemma~2.1.}~~For $x,y \in S$, write $x \sim y$ if
$x$ and $y$ lie in exactly the same collection of translates of $X$.
This clearly defines an equivalence relation on $S$, and our goal,
of course, is to show that $u \not\sim v$

Since $X$ is nonempty and proper in $S$, there exist inequivalent
members of $S$, and so the equivalence class $\Delta$ of $u$ is
proper in $S$. The action of $G$ permutes the equivalence classes,
and it is clear that $\Delta$ is a block. By primitivity, therefore,
$\Delta = \{u\}$, and thus $u \not\sim v$, as wanted.\qed
\medskip

The conclusion of Lemma~2.1 is only slightly weaker than the full
conclusion of Rudio's lemma, and we digress briefly to explain this.
Suppose we know that some translate of $X$ contains exactly one
of $u$ or $v$ and we wish to be able to specify that (say) $u$ is in a
translate of $X$ and $v$ is not. We show that this is possible if we
assume one additional piece of information: that $ug = v$ for some
element $g \in G$ of {\it finite} order. Of course, this condition
holds for finite groups $G$, and thus together with Lemma~2.1, the
following argument provides an alternative proof of Rudio's lemma.

Suppose that some translate $Y$ of $X$ contains $v$ but not $u$
and assume that $v = ug$, where $g^n = 1$ for some positive
integer $n$. We show that an appropriate translate of $Y$ contains
$u$ but not $v$. We have $ug = v \in Y$ but $ug^n = u \not\in Y$,
and it follows that there exists an integer $m$ such that
$ug^m \in Y$ but $ug^{m+1} \not\in Y$. Then $u \in Yg^{-m}$ but
$v \not\in Yg^{-m}$ because $vg^m = ug^{m+1} \not\in Y$.

We return now to our main theme and establish the inequality
$(t+1)d \ge 2n$ when $1 \le t < n$ and $S$ is a primitive finite
$G$-set. 
\medskip

\iitem{Proof of Theorem~C.}~~We use a refinement of the
argument in the proof of Theorem~B. Choose a basis $\BB$ for
$\gen v$ that consists of translates of $v$ and note that
$\bigcup_{b \in \BB} \supp b = S$. For each point $s \in S$, write
$m(s)$ to denote the number of members $b \in \BB$ such that
$s \in \supp b$, and observe that $m(s) \ge 1$ for all $s \in S$. Since
$|\BB| = d$ and $|\supp b| = t$ for all $b \in \BB$, we see that
$td = \sum_{s \in S} m(s)$. (This is because each side of this
equation counts the number of ordered pairs
$(b,s) \in \BB \times S$ such that $s \in \supp b$.)

We claim that for each member $b \in \BB$, there is at most one
point $s \in \supp b$ such that $m(s) = 1$. Assuming this for the
moment, we see that there are at most $|\BB|$ points $s \in S$
such that $m(s) = 1$ and that for all other points $s \in S$, we have
$m(s) \ge 2$. Then
$$
td = \sum_{s \in S} m(s) \ge 2|S| - |\BB| = 2n - d \leqno{(*)}
$$
and the desired inequality follows.

Suppose then, that there exists $b \in \BB$ and distinct points
$x,y \in \supp b$ such that $m(x) = 1 = m(y)$. Then $b$ is the only
member of $\BB$ for which either $x$ or $y$ is in the support. It
follows that for each linear combination $c$ of vectors in $\BB$,
either the support of $c$ contains both $x$ and $y$ or neither of
them, depending on whether or not $b$ appears with nonzero
coefficient in the expansion of $c$ in terms of $\BB$.

Now if $w$ is any $G$-translate of $v$, we have $w \in \gen v$, and
thus $w$ is a linear combination of vectors in $\BB$. It is not the
case, therefore, that $\supp{w}$ contains exactly one of $x$ and
$y$. Writing $X = \supp v$, it follows that no translate of $X$
contains exactly one of $x$ and $y$. But $1 \le t < n$, and thus $X$
is nonempty and proper in $S$, and since $S$ is primitive, this
contradicts Lemma~2.1. The proof is now complete.\qed
\medskip

In fact, this argument gives some additional information, which we
will exploit in the following section.
\medskip

\iitem{(2.2) LEMMA.}~~{\sl Suppose that equality holds in
Theorem~C and let $\BB$ be any linearly independent set of
translates of $v$. Then each point of $S$ is in the support of at
most two members of $\BB$. Also, the support of each member of
$\BB$ contains a point that is not in the support of any other
member of this set.}
\medskip

\iitem{Proof.}~~Since $\BB$ is part of some basis for $\gen v$
consisting of translates of $v$, it is no loss to assume that $\BB$ is
such a basis. Since we are assuming that $(t+1)d = 2n$, it follows
that equality holds in $(*)$, in the proof of Theorem~C. But we
know that $m(s) \ge 1$ for all $s \in S$ and that $m(s) = 1$ for at
most one point in the support of each of the $d$ members of
$\BB$. We see, therefore, that equality forces $m(s) = 1$ for
exactly one point in the support of each member of $\BB$ and
$m(s) = 2$ for all other points. This completes the proof.\qed
\bigskip

\iitem{3. Equality in Theorem~C.}

When is it true that the equality $(t+1)d = 2n$ holds for some
vector $v \in F[S]$, where $S$ is a primitive $G$-set? If $t = 1$,
then by Theorem~B, we have $d = n$, and equality automatically
holds. At the opposite extreme, when $t = n$, it is clear that
equality is impossible except in the degenerate case where $n = 1$.
We will discuss the case where $t = n-1$ later, and so we assume
for this section that $1 < t < n-1$, and we work toward a proof of
Theorem~D, which gives highly restrictive necessary conditions for
equality to hold. As we shall see, these conditions turn out to be
sufficient too.
\medskip

\iitem{(3.1) LEMMA.}~~{\sl Suppose that equality holds in
Theorem~C and let $X$ and $Y$ be distinct translates of $\supp v$.
Then the stabilizer in $G$ of each point of $X \cap Y$ stabilizes the
set $X \cup Y$.}
\medskip

\iitem{Proof.}~~Let $a$ and $b$ be translates of $v$ that are
supported on $X$ and $Y$, respectively, and note that $a$ and $b$
are linearly independent since they have distinct supports. Suppose
that $x \in X \cap Y$ and that $xg = x$, where $g \in G$. Then $x$ is
in the support of $a$, $b$ and $ag$, and hence by Lemma~2.2,
these three vectors cannot be linearly independent. We conclude
that $ag$ is a linear combination of $a$ and $b$, and so
$\supp{ag} \sbs \supp a \cup \supp b$. Thus $Xg \sbs X \cup Y$, and
similarly, $Yg \sbs X \cup Y$. It follows that
$(X \cup Y)g = X \cup Y$, as required.\qed
\medskip

We also need the following elementary combinatorial fact and an
easy consequence. 
\medskip

\iitem{(3.2) LEMMA.}~~{\sl Let $\Omega$ be a collection of subsets
of some set $S$ and assume that each point of $S$ lies in exactly
two members of $\Omega$. Say that a subcollection
$\Lambda \sbs \Omega$ is {\bf even} if every point of $S$ lies in an
even number of members of $\Lambda$. The following then hold.
\smallskip
\ritem{(a)} The intersection of any two even subcollections of
$\Omega$ is even.
\smallskip
\ritem{(b)} If $\Lambda \sbs \Omega$ is even, then each member of
$\Lambda$ is disjoint from each member $\Omega - \Lambda$.
\medskip}

\iitem{Proof.}~~Suppose that $\Lambda, \Delta \sbs \Omega$ are
even. Let $s \in S$ and suppose that the number of members of
$\Lambda \cap \Delta$ that contain $s$ is nonzero. If $X$ and $Y$
are the two members of $\Omega$ that contain $s$. Then at least
one of $X$ or $Y$ is in both $\Lambda$ and $\Delta$, and thus
since these collections are even, both $X$ and $Y$ lie in both
$\Lambda$ and $\Delta$. Then $s$ lies in exactly two members of
$\Lambda \cap \Delta$, establishing (a).

Now suppose that $s \in S$ lies in some member $X \in \Lambda$,
where $\Lambda$ is even. Then both members of $\Omega$ that
contain $s$ lie in $\Lambda$, and so no member of
$\Omega - \Lambda$ can contain $s$. This proves (b).\qed
\medskip

\iitem{(3.3) LEMMA.}~~{\sl Let $S$ be a finite primitive $G$-set and
let $\TT$ be a set of vectors in $F[S]$, where $F$ is an arbitrary
field. Suppose that the members of $\TT$ have distinct supports
and that $\Omega = \{\supp v \mid v \in \TT\}$ is transitively
permuted by $G$. If every point of $S$ lies in exactly two members
of $\Omega$, then every proper subset of $\TT$ is linearly
independent, and so the dimension of the linear span of of $\TT$ is
either $|\TT|$ or $|\TT| - 1$.}
\medskip

\iitem{Proof.}~~We can suppose that $\TT$ is not linearly
independent, and we choose a minimal dependent subset
$\TT_0 \sbs \TT$. We can thus write $\sum a_v v = 0$, where the
sum runs over $v \in \TT_0$ and all of the coefficients $a_v \in F$
are nonzero. It follows that no point of $S$ can lie in the support
of exactly one vector $v \in \TT_0$.

Let $\Omega_0 = \{\supp v \mid v \in \TT_0\}$. We have seen that
each point of $S$ lies in either zero or two members of
$\Omega_0$, and so the language of Lemma~3.2, the collection
$\Omega_0$ is even. Since $\Omega_0$ is nonempty, we can choose
a minimal nonempty even subcollection $\Lambda \sbs \Omega_0$. 

Let $U = \bigcup\Lambda$, the union of all members of $\Lambda$.
Suppose that $g \in G$ is an element such that $Ug \ne U$, and
observe that $\Lambda g \ne \Lambda$. Now $\Lambda g$ is even,
and thus $\Lambda \cap \Lambda g$ is even by Lemma~3.2(a). But
$\Lambda \cap \Lambda g$ is proper in $\Lambda$ and it follows by
the minimality of $\Lambda$, that
$\Lambda \cap \Lambda g = \emptyset$. Thus no member of
$\Lambda g$ is in $\Lambda$, and so by Lemma~3.2(b), we see that
each member of $\Lambda g$ is disjoint from $U$. It follows that
$U \cap Ug = \emptyset$, and hence $U$ is a block.

Since each point of $S$ is in more than one member of $\Omega$,
the members of $\Omega$ cannot be singleton sets, and hence the
block $U$ is not a singleton. But $G$ is primitive on $S$, and so we
must have $U = S$, and thus no member of $\Omega$ is disjoint
from $U$. By Lemma~3.2(b) it follows that $\Lambda = \Omega$,
and thus $\Omega_0 = \Omega$. Then $\TT_0 = \TT$, and hence
every proper subset of $\TT$ is linearly independent. The span of
$\TT$, therefore, has dimension $|\TT| - 1$, as wanted.\qed
\medskip

Recall that Theorem~D asserts (among other things) that the group
$G$ has a doubly transitive action in which the stabilizer of a
two-point subset is a maximal subgroup. (In other words, its action
on the two-point subsets is primitive.) In order to obtain the
conclusion in Theorem~D that the field $F$ must have
characteristic $2$, we need to study such doubly transitive groups,
and that is the purpose of the following lemma. This lemma will also
be crucial in finding the examples where equality holds in
Theorem~C.
\medskip

\iitem{(3.4) LEMMA.}~~{\sl Let $G$ be a doubly transitive
permutation group on a set $\Omega$ and assume that the stabilizer
of a two-point subset of $\Omega$ is a maximal subgroup of $G$. If
$N$ is a minimal normal subgroup of $G$, then either $N$ is a
nonabelian simple group or $|N| = 3 = |\Omega|$.}
\medskip

The key to the proof of Lemma~3.4 is to show that $N$ must be
primitive on $\Omega$, and then it follows by a standard argument
that $N$ is either abelian or simple. In the abelian case, the
conclusion that $|N| = 3$ follows easily from the primitivity of $G$
on the two-point subsets of $\Omega$. It is true (but not quite
trivial) that in general, a nonabelian minimal normal subgroup of a
$2$-transitive group must be primitive, and hence simple. (See
page 202 of \ref\bur.) The primitivity of the action of $N$ on
$\Omega$ is much easier to prove, however, in our case, where
$G$ acts primitively on the two-point subsets. We have decided,
therefore, to give the direct and elementary proof. 
\medskip

\iitem{Proof of Lemma~3.4.}~~Write $|\Omega| = r$ and observe
that $r \ge 3$ since the stabilizer of a two-point subset is proper.
Also, because $G$ is $2$-transitive on $\Omega$, we see that $N$ is
transitive, and thus $N$ cannot fix a two-point subset of $\Omega$.
Since the action of $G$ on these subsets is is primitive and $N$ acts
nontrivially, it follows that $N$ is transitive on the collection of
two-point subsets, and this implies that $N$ is primitive on
$\Omega$. (To see this, observe that if there were a nontrivial
$N$-block $\Delta \sbs \Omega$, we could choose
$\alpha,\beta \in \Delta$ and $\gamma \in \Omega - \Delta$, and
then no element of $N$ could take $\{\alpha,\beta\}$ to
$\{\alpha,\gamma\}$.)

Since $N$ is minimal normal in $G$, we see that if it is not simple,
we can write $N = A \times B$, where each of $A$ and $B$ is
nontrivial. Since $N$ is primitive on $\Omega$, it follows that each
of $A$ and $B$ is transitive, and since $A$ and $B$ centralize each
other, we conclude that they are both regular. Thus
$|A| = r = |B|$ and $|N| = r^2$.

Now let $K = G_\alpha$ be a point stabilizer in $G$ and note that
$|K \cap N|$ is coprime to $r - 1$ since $|N| = r^2$. Also, $K$ acts
transitively on the $r - 1$ points of $\Omega - \{\alpha\}$, and
$K \cap N \nor K$. It follows from this that $K \cap N$ acts trivially
on $\Omega$. Then $K \cap N = 1$, and so $N$ is regular and
$|N| = r$, which is a contradiction. We conclude that $N$ is simple.

Finally, if $N$ is abelian, then it is regular and $|N| = r$. But $N$
acts transitively on the $r(r - 1)/2$ two-point subsets of $\Omega$,
and thus $r \ge r(r-1)/2$ and $r \le 3$. We know that $r \ge 3$,
however, and thus $r = 3$, as required.\qed  
\medskip

\iitem{(3.5) COROLLARY.}~~{\sl Let $G$ be a doubly transitive
permutation group on a set $\Omega$ having at least four points and
assume that the stabilizer of a two-point subset of $\Omega$ is a
maximal subgroup of $G$. If $1 < N \nor G$, then $N$ acts doubly
transitively on $\Omega$.}
\medskip

\iitem{Proof.}~~As in the previous proof, $N$ acts nontrivially on the
collection of two-point subsets of $\Omega$. Since $G$ is primitive
on this collection and $N$ is a normal subgroup that acts nontrivially,
it follows that $N$ acts transitively on the collection of two-point
subsets. To show that $N$ acts doubly transitively on $\Omega$,
therefore, it suffices to show that some element of $N$ interchanges
some pair of points of $\Omega$. By Lemma~3.4, however, $N$ is not
solvable, and thus by the Feit-Thompson theorem, $N$ contains an
element of order $2$. Such an element, of course, interchanges a
pair of points.\qed
\medskip

The following result is a somewhat more precise version of
Theorem~D of the introduction.
\medskip

\iitem{(3.6) THEOREM.}~~{\sl Suppose that equality holds in
Theorem~C and that $1 < t < n-1$. Let $\Omega$ be the set of
$G$-translates of $\supp v$. The following then hold.
\smallskip
\ritem{(a)} Each point of $S$ lies in exactly two members of
$\Omega$ and every two distinct members of $\Omega$ have
exactly one point of $S$ in common. 
\smallskip
\ritem{(b)} The action of $G$ on $\Omega$ is $2$-transitive and the
induced action on two-element subsets of $\Omega$ is primitive.
\smallskip
\ritem{(c)} $d = |\Omega| - 1 = t$.
\smallskip
\ritem{(d)} Each member of $\Omega$ is transitively permuted by
its (setwise) stabilizer in $G$.
\smallskip
\ritem{(e)} $v$ is a scalar multiple of the sum of the points in its
support.
\ritem{(f)} $F$ has characteristic $2$.}
\medskip

\iitem{Proof.}~~We have $(t+1)d = 2n$ and $2 < t + 1 < n$, and
thus $2 < d < n$. In particular, since $d$ is a divisor of $2n$, we see
that $n$ cannot be prime. Also, since $G$ is primitive on $S$ and
$|S|$ is not prime, it is easy to see that $G$ is generated by the
stabilizers of any two distinct points of $S$.

Let $X,Y \in \Omega$ be distinct and write $u = |X \cap Y|$. We
argue first that $u \le 1$. Otherwise, let $x,y \in X \cap Y$ be
distinct. By Lemma~3.1, the point stabilizers $G_x$ and $G_y$ both
stabilize the set $X \cup Y$, and since these subgroups generate
$G$, it follows that $G$ stabilizes $X \cup Y$. Then $X \cup Y = S$,
and we have $n = 2t - u$.

Now let $Z,W \in \Omega$ be arbitrary. Then
$Z \cup W \sbs S = X \cup Y$, and since all members of $\Omega$
have equal cardinality, it follows that
$|Z \cap W| \ge |X \cap Y| \ge 2$. The previous argument now
shows that if $Z \ne W$, then $S = Z \cup W$, and thus
$|Z \cap W| = |X \cap Y| = u$.

Since the set $X \cap Y$ is not $G$-invariant, it cannot be the case
that $\Omega = \{X,Y\}$, and thus we can choose $Z \in \Omega$,
different from $X$ and $Y$. Then $X \cup Z = S = X \cup Y$, and so
$S - X \sbs Y \cap Z$. Then $n - t = |S - X| \le |Y \cap Z| = u$, and
we have $2t - u = n \le t + u$, and thus $t \le 2u$. We have
$$
td < (t+1)d = 2n = 4t - 2u \le 3t \,,
$$
and we conclude that $d < 3$. Since $d > 2$, this is a contradiction,
and thus every two distinct members of $\Omega$ have at most one
point in common.

Now let $X,Y,Z \in \Omega$ be distinct and suppose that
$X \cap Y \cap Z$ is nonempty. Let $x \in X \cap Y \cap Z$ and
observe that by the result of the previous argument, the intersection
of any two of $X$, $Y$ and $Z$ is exactly the set $\{x\}$. Let $a$, $b$
and $c$ be translates of $v$ with supports $X$, $Y$ and $Z$,
respectively, and note that by Lemma~2.2, these three vectors
cannot be linearly independent. Thus $c$ (say) is a linear
combination of $a$ and $b$, and hence $Z \sbs X \cup Y$. Then
$Z - X \sbs Y$, and we have $Z - X \sbs Y \cap Z = \{x\}$. But $x \in X$,
and it follows that $Z - X$ is empty. This is a contradiction since
$|X| = |Z|$ and $X \ne Z$. We conclude that each point of $S$ lies in
at most two members of $\Omega$. 

The action of $G$ on $S$ is primitive, and since $1 < t < n$, it
follows that $\supp v$ is not a block. The members of $\Omega$,
therefore, are not pairwise disjoint, and so some point of $S$ lies in
two of them. Because $G$ is transitive on $S$, it follows that every
point in $S$ lies in exactly two members of $\Omega$, and this
defines a map from $S$ into the collection of two-element subsets
of $\Omega$. This map is injective since distinct members of
$\Omega$ can have at most one point in common. To complete the
proof of (a), we must show that our map is surjective, and so it
suffices to show that $r(r-1)/2 = n$, where $r = |\Omega|$. Since
$2n = (t+1)d$, we see that (a) will follow once we prove (c).

Since $|\Omega| = r$, we see that $2n = tr$ because both sides of
this equation count pairs $(x,X) \in S \times \Omega$ such that
$x \in X$. Then $tr = (t+1)d$, and so $d < r$. If two translates of $v$
have the same support, then by Lemma~2.2, they must be scalar
multiples of one another. It follows that $\gen v$ is spanned by a
collection of translates of $v$, one supported on each of the $r$
members of $\Omega$. Since $d < r$, we see by Lemma~3.3 that
$d = r - 1$, as wanted. Then $(t+1)d = 2n = tr = t(d+1)$, and it
follows that $t = d$. This establishes (c) and (a).

We have a natural bijection between $S$ and the collection of
two-member subsets of $\Omega$, and so we know that $G$ is
transitive, and in fact primitive, on this collection. To show
that $G$ is $2$-transitive on $\Omega$ and thus complete the proof
of (b), it suffices to show that if $x \in S$, then some element of
$G_x$ interchanges the two members $X,Y \in \Omega$ that contain
$x$. Certainly, $G_x$ permutes $\{X,Y\}$, and it suffices,
therefore, to show that $G_x$ does not stabilize $X$. But if
$G_x \sbs G_X$, then since $G_x$ is maximal in $G$, we have
$G_x = G_X$ and $tr/2 = n = |G:G_x| = |G:G_X| = r$. Then $t = 2$,
and since we know that $t = d \ne 2$, we have a contradiction. This
proves (b).

Let $X \in \Omega$. To prove (d), we must show that $G_X$ is
transitive on $X$. But the points of $X$ are exactly the
intersections of $X$ with the various members of $\Omega - \{X\}$,
and since $G$ is $2$-transitive on $\Omega$, these sets are
transitively permuted by $G_X$, and thus the points of $X$ are also
transitively permuted by $G_X$, as wanted.

For (e), let $X = \supp v \in \Omega$ and write
$K = G_X$. We know by Lemma~2.2 that $vk$ must be a scalar
multiple of $v$ for each element $k \in K$, and we define the map
$\lambda:K \to F$ by $vk = \lambda(k)v$. Observe that $\lambda$ is
a homomorphism from $K$ into $F^\times$.

Let $A$ be a set of representatives for the right cosets of $K$ in
$G$ and note that the $r = d+1$ vectors $va$ for $a \in A$ have
distinct supports. Also, every translate of $v$ is a scalar multiple of
one of the vectors $va$ with $a \in A$, and hence these vectors
span $\gen v$, which has dimension $d$. It follows that there is a
unique (up to scalar multiplication) linear dependence relation
among the vectors $va$, and we write $\sum_{a \in A} c_a va = 0$,
where the coefficients $c_a$ lie in $F$. Also, we know from
Lemma~3.3 that every proper subset of $\{va \mid a \in A\}$ is
linearly independent, and thus $c_a \ne 0$ for all $a \in A$.

Let $g \in G$ be arbitrary and let $a \in A$. Write $a \mdot g$ to
denote the unique element of $A$ that lies in the coset $Kag$. Now
apply $g$ to the equation $\sum c_a va = 0$ and express each
vector $vag$ as an appropriate scalar multiple of $v(a \mdot g)$. (In
fact, it is easy to see that this scalar is $\lambda(ag(a \mdot g)\inv)$,
although we shall not need this explicit formula.) Since the
map $a \mapsto a \mdot g$ is a permutation of $A$, what results is
a new dependence relation of the form $\sum b_a va = 0$. It follows
that there must exist a scalar $\mu(g)$, depending only on $g$,
such that $b_a = \mu(g) c_a$ for all $a \in A$. Furthermore, it is not
hard to see that $\mu: G \to F^{\times}$ is a group homomorphism.
(What is really going on here is that $Fv$ is a $1$-dimensional
$K$-module, and the vector $\sum_{a \in A} c_a (v \tens a)$ spans a
$1$-dimensional $G$-submodule of the induced module $(Fv)^G$.) 

We can suppose that $1 \in A$. If we take $g \in K$, then
$1 \mdot g = 1$, and since $vg = \lambda(g) v$, we see that
$b_1 = \lambda(g) c_1$, and thus $\mu(g) = \lambda(g)$ and $\mu$
is an extension of $\lambda$ to $G$. (This conclusion could also be
proved using an appropriate generalization of Frobenius
reciprocity.)

Write $N = \ker\mu \nor G$ and observe that $N > 1$ since it is clear
that $G$ is noncyclic, and thus $N$ is doubly transitive on $\Omega$
by Corollary~3.5. Also, $N \cap K = \ker\lambda$, and this is exactly
the stabilizer of $v$ in $K$, and hence it is the stabilizer of $v$ in
$G$.

Let $x \in X$ and write $H = G_x$. If $g \in H \cap K$, then
$vg = \lambda(g)v$ and the (nonzero) coefficients of $x$ in $v$ and
$vg$ are equal. It follows that $\lambda(g) = 1$, and thus
$\mu(g) = 1$ and $g \in N$. This shows that $H \cap K \sbs N$.

Now $H \cap K$ stabilizes both $X$ and also the unique other
member of $\Omega$ that contains $x$. In fact, $H \cap K$ is a full
two-point stabilizer in the $2$-transitive action of $G$ on
$\Omega$, and so $|G:H \cap K| = r(r-1)$. If $N < G$, then
$|N:H \cap K| < r(r-1)$, and thus $N$ is not $2$-transitive on
$\Omega$, and this is a contradiction. We conclude, therefore,
that $N = G$, and thus $K = N \cap K$ stabilizes $v$. By (d),
however, $K$ is transitive on $X = \supp v$, and (e) follows.

We can replace $v$ by a scalar multiple and assume that $v$ is
exactly the sum of the points in its support. For each member
$X \in \Omega$, write $v_X = \sum_{x \in X} x$, and note that these
$r = |\Omega|$ vectors are exactly the $G$-translates of $v$. But
$d < r$, and so the vectors $v_X$ are dependent, and we can write
$\sum_{X \in \Omega} c_X v_X = 0$ for suitable coefficients
$c_X \in F$, not all $0$.

Choose distinct members $X,Y,Z \in \Omega$ with $c_X \ne 0$.
There is a point of $S$ that lies in $X$ and $Y$ and in no other
member of $\Omega$, and it follows that $c_X + c_Y = 0$. Similarly,
$c_X + c_Z = 0$ and $c_Y + c_Z = 0$, and it follows that $2c_X = 0$.
We deduce that $F$ has characteristic $2$, as required for (f). The
proof is now complete.\qed
\medskip

To summarize, we see that if equality holds in Theorem~C and
$1 < t < n-1$, then $G$ is a $2$-transitive group on some set
$\Omega$ such that the induced action on the two-point subsets of
$\Omega$ is primitive. Also, $t = d$, and thus $2n = t(t + 1)$. It
follows that we cannot have $t = 2$ since otherwise $n = 3$ and the
inequality $t < n - 1$ would not hold. Thus $t \ge 3$, and so
$|\Omega| = t + 1 \ge 4$.

It is possible (by appealing to the classification of finite simple
groups) to list all possible doubly transitive groups $G$ on four or
more points such that the action of $G$ on the two-point subsets is
primitive. We shall see in the next section that given any such group
and any field of characteristic $2$, it is possible to construct a
corresponding example where equality holds in Theorem~C.
\bigskip

\iitem{4. The examples where $1 < t < n-1$.}

We know that if equality holds in Theorem~C and $1 < t < n-1$,
then $G$ has a doubly transitive permutation representation on
some set $\Omega$ of cardinality $t + 1$, where the induced
action on two-point subsets of $\Omega$ is primitive.

Conversely, suppose that $G$ is doubly transitive on some set
$\Omega$ with $|\Omega| \ge 3$. Let $S$ be the set of two-point
subsets of $\Omega$ and suppose that the action of $G$ on $S$ is
primitive. Let $K$ be the stabilizer in $G$ of a point
$\alpha \in \Omega$ and let $H$ be the stabilizer in $G$ of the
two-point set $\{\alpha,\beta\} \in S$. 

Since $G$ is doubly transitive on $\Omega$, some element $t \in G$
interchanges $\alpha$ and $\beta$, and we see that
$t \in H$, and so $H$ is transitive on $\{\alpha,\beta\}$. It follows
that $|H:H \cap K| = 2$. Also, $K$ is transitive on
$\Omega - \{\alpha\}$, which contains at least $2$ points, and so
$K$ does not fix $\beta$ and it follows that $K \not\sbs H$. Finally,
we observe that since $|\Omega| > 2$, we have $|S| > 1$, and so
$K$, which stabilizes a point of $S$ is certainly not transitive on
$S$. 

We consider a somewhat more general situation.
\medskip

\iitem{(4.1) LEMMA.}~~{\sl Let $S$ be a finite primitive right
$G$-set and let $H = G_x$ for some point $x \in S$. Suppose that
$K \sbs G$ is any subgroup that satisfies the following three
conditions:
\smallskip
\ritem{(1)} $|H:H \cap K| = 2$.
\smallskip
\ritem{(2)} $K \not\sbs H$.
\smallskip
\ritem{(3)} $K$ is not transitive on $S$.
\smallskip
Let $F$ be any field, and let $v \in F[S]$ be the sum of the points in
the $K$-orbit of $x$. As usual, write $t = t(v)$, $d = d(v)$ and
$n = |S|$. The following then hold.
\smallskip
\ritem{(a)} $t = |K:H \cap K|$.
\smallskip
\ritem{(b)} $t|G:K| = 2n$.
\smallskip
\ritem{(c)} Either $d = |G:K|$ or $d = |G:K| - 1$.
\smallskip
\ritem{(d)} If $F$ has characteristic $2$, then $d = |G:K| - 1$.}
\medskip

\iitem{Proof.}~~Statement (a) is clear since $|K:H \cap K|$ is the
size of the $K$-orbit of $x$. Thus
$t|G:K| = |G:H \cap K| = 2|G:H| = 2n$, where the second equality
follows by assumption (1), and this establishes (b). We work next to
determine $d$.

We claim that $K$ is the full stabilizer in $G$ of the $K$-orbit $X$
containing $x$. Otherwise, there exists a subgroup $J \sbs G$ such
that $J$ stabilizes $X$ and $K < J$. In fact, $J < G$ since $X$ is
proper in $S$ by assumption (3). Since $K$ is transitive on $X$ we
have $J = KJ_x = K(H \cap J)$, and thus $H \cap J \not\sbs K$. Then
$H \cap K < H \cap J \sbs H$. But $|H:H \cap K| = 2$, and so it
follows that $H \cap J = H$ and $H \sbs J$. In fact, $H < J$ since
$K \sbs J$ but $K \not\sbs H$ by assumption (2). We now have
$H < J < G$, and this is a contradiction since by the primitivity of
$S$, the point stabilizer $H$ is maximal in $G$.

Let $\Omega$ be the set of $G$-translates of $X$. Then
$|\Omega| = |G:K|$ since we now know that $G_X = K$. Since $G$
is transitive on $S$ and on $\Omega$, we see that the number of
members of $\Omega$ that contain each point $y \in S$ is some
constant $m$, independent of the choice of $y$. If we count
ordered pairs $(y, Y)$, where $Y \in \Omega$ and $y \in Y$, we see
that $mn = t|\Omega| = t|G:K|$, and thus $m = 2$ by (b). Each point
of $S$, therefore, lies in exactly two members of $\Omega$.

The translates of $v$ are exactly the $|G:K|$ vectors
$v_Y = \sum_{y \in Y} y$ for $Y \in \Omega$, and $d$ is the
dimension of the linear span of these vectors. By Lemma~3.3,
therefore, either $d = |G:K|$ or $d = |G:K| - 1$, proving (c).

 If $F$ has characteristic $2$, it follows from the fact that each
point of $S$ lies in exactly two members of $\Omega$ that
$\sum v_Y = 0$, and thus the vectors $v_Y$ are dependent and we
have $d < |G:K|$. Then $d = |G:K| - 1$, as wanted.\qed 
\medskip

In the situation of Lemma~4.1, write $r = |G:K|$. Then
$d = r$ or $d = r - 1$ and $tr = 2n$. Thus either
$td = 2n$ or $t(d+1) = 2n$, and if $F$ has characteristic $2$, then
only the second alternative can occur.

In the case where $G$ is doubly transitive on $|\Omega|$ and $H$
is the stabilizer of a two-point subset, we have $r = |\Omega|$, and
thus $n = (r-1)r/2$. If $F$ has characteristic $2$, then
$d = r - 1$ and $t(d+1) = 2n = (r - 1)r = d(d + 1)$, and so $t = d$. Then
$(t + 1)d = 2n$ and we have equality in Theorem~C. Note also that if
$|\Omega| \ge 4$, then $t = r - 1 \ge 3$, and we certainly have
$1 < t < n-1$.

We now address the question of finding all doubly transitive
permutation groups $G$ acting on a set $\Omega$ consisting of at
least four points and such that the induced action of $G$ on the
two-point subsets of $\Omega$ is primitive.

First, recall Lemma~3.4, which asserts that a minimal normal
subgroup $N$ of $G$ must be a nonabelian simple group. By
Corollary~3.5, furthermore, the simple group $N$ must itself be
doubly transitive on $\Omega$, and so $G$ is contained between
$N$ and its normalizer in the symmetric group $S_\Omega$. (Note
that since $N$ clearly has trivial centralizer in $S_\Omega$, its
normalizer, which we denote ${\rm Aut}_\Omega(N)$, is naturally
embedded in $\aut N$, and so can be computed.)

To find the examples we seek, therefore, we start by examining
the list of doubly transitive simple groups, as compiled by
P.~Cameron \ref\cam. For each such simple group $N$, acting
doubly transitively on a set $\Omega$, we compute
$A = {\rm Aut}_\Omega(N)$ and we check the groups $G$ with
$N \sbs G \sbs A$ to determine which of them (if any) act
primitively on the two-point subsets of $\Omega$. Once we find
such a group $G$, then of course, any larger subgroup (contained
in $A$) will also yield an example.

The following table lists all simple groups $N$ that occur as minimal
normal subgroups of doubly transitive groups that are primitive on
the set of two-point subsets. In most of these cases, the simple
group $N$ is itself primitive on the two-point subsets, and so every
group contained between $N$ and ${\rm Aut}_\Omega(N)$ is an
example, as desired. The only exceptions are $PSL(2,q)$ acting on
$q + 1$ points, where $q \in \{7,9,11\}$. These three simple groups
do not act primitively on the two-point subsets, but the overgroup
$PGL(2,q)$ does have this property. (And for $q = 9$, the
overgroup $M_{10}$ provides another example.)
\medskip

\centerline{\vbox{\offinterlineskip
\halign{\strut
\vrule\hfil~#~\hfil&\vrule\hfil~#~\hfil&\vrule\hfil~#~\hfil\vrule\cr
\noalign{\hrule}
Simple group&Degree&Remarks\cr
\noalign{\hrule\vglue1truept\hrule}
$A_r$ (Alternating)&$r$&$r \ge 5$\cr
\noalign{\hrule}
$M_r$ (Mathieu)&$r$&$r \in \{11,12,22,23,24\}$\cr
\noalign{\hrule}
$PSL(2,q)$&$q + 1$&$q \not\in \{2,3,5\}$\cr
\noalign{\hrule}
$Sz(q)$ (Suzuki)&$q^2 + 1$&$q = 2^{2k+1}~~~k \ge 1$\cr
\noalign{\hrule}
$PSL(2,11)$&$11$&\cr
\noalign{\hrule}
$M_{11}$ (Mathieu)&$12$&\cr
\noalign{\hrule}
$HS$ (Higman-Sims)&176&\cr
\noalign{\hrule}
$Co_3$ (Conway)&276&\cr
\noalign{\hrule}
}}}
\medskip

We will not actually prove that this table is complete, and we give
only a partial proof that it is correct. First, (and this does not
depend on the classification of simple groups) almost every group
that is triply (and not just doubly) transitive and that fails to have
an abelian normal subgroup is guaranteed to act primitively on
two-point subsets. (The exception here is the action of the
symmetric group $S_5$ on $6$ points.) This justifies the appearance
in our table of the alternating groups and the Mathieu groups in
their natural permutation representations; it covers the groups
$PSL(2,q)$, where $q$ is a power of $2$ exceeding $2$, and also it
shows that $M_{11}$ in its $12$-point representation belongs in our
table.

This fact about triply transitive groups is an easy corollary of the
following theorem of Cameron \ref\camtt.
\medskip

\iitem{(4.2) THEOREM (Cameron).}~~{\sl Suppose that $G$ is a triply
transitive permutation group on a set $\Omega$ and let
$\alpha \in \Omega$. Assume that $G$ also acts on a set $\Lambda$
in such a way that $\Omega - \{\alpha\}$ and $\Lambda$ are
isomorphic as $G_\alpha$-sets. Then either $G$ has a nontrivial
abelian normal subgroup, or else $G \cong S_5$ and
$|\Omega| = 6$.}
\medskip

\iitem{(4.3) COROLLARY.}~~{\sl Let $G$ be a triply transitive
permutation group on a set $\Omega$ and assume that $G$ has no
nontrivial abelian normal subgroup. Then the action of $G$ on the
two-point subsets of $\Omega$ is primitive unless $G \cong S_5$ and
$|\Omega| = 6$.}
\medskip

\iitem{Proof.}~~Let $\alpha,\beta \in \Omega$ and (as usual) let
$K = G_\alpha$ and $H = G_{\{\alpha,\beta\}}$. Assuming that the
action of $G$ on the two-point subsets of $\Omega$ is not
primitive, the subgroup $H$ is not maximal, and we let $H < J < G$.
Since $G$ is triply transitive, we see that $\{\alpha,\beta\}$ and
$\Omega - \{\alpha,\beta\}$ are the orbits of $H$ on $\Omega$. But
$J$ does not stabilize the set $\{\alpha,\beta\}$ and $J > H$, and
it follows that $J$ is transitive on $\Omega$ and we have $JK = G$. 

Since $J < G$, we see that $K \not\sbs J$, and thus
$H \cap K \sbs J \cap K < K$. But $K$ acts doubly transitively, and
hence primitively on $\Omega - \alpha$ and $H \cap K$ is the
stabilizer of the point $\beta$ in this action. It follows that
$H \cap K$ is maximal in $K$ and we conclude that
$H \cap K = J \cap K$.

Let $\Lambda$ be the set of right cosets of $J$ in $G$. Then $K$
acts transitively on $\Lambda$ since $JK = G$, and
$K \cap J = K \cap H$ is the stabilizer of a point in this action. Also,
$K$ acts transitively on $\Omega - \{\alpha\}$ and $K \cap H$ is the
stabilizer of a point in that action. It follows that $\Lambda$ and
$\Omega - \{\alpha\}$ are isomorphic as $K$-sets, and the result
follows via Cameron's theorem.\qed
\medskip

Finally, we explain why the group $PSL(2,q)$ appears in our table
when $q > 11$ is odd. Consider the action of $SL(2,q)$ on the
$q + 1$ subspaces of dimension $1$ in a $2$-dimensional space $V$
over the field $F$ of order $q$. If $\alpha = Fv$ and
$\beta = Fw$ are distinct ``points" and we use $\{v,w\}$ as our basis
for $V$, then the stabilizer in $SL(2,q)$ of the set
$\{\alpha,\beta\}$ is easily seen to be the group of monomial
matrices of determinant $1$. This group is dihedral of order
$2(q - 1)$, and it corresponds to a dihedral subgroup of order
$q - 1$ in $PSL(2,q)$. Since $q > 11$, we see that neither $A_5$ nor
$S_4$ contains a dihedral subgroup of order $q - 1$, and it follows
by checking the list of isomorphism types of subgroups of
$PSL(2,q)$ that a dihedral subgroup of order $q - 1$ is necessarily
maximal. (See \ref\hup, Hauptsatz II.8.27.) In other words, the
action of $PSL(2,q)$ on the two-point subsets of the projective line
over $F$ is primitive, as claimed.

We close this section with another application of Lemma~4.1. We
know by Theorem~C that $(t+1)d \ge 2n$ in the primitive case, and
we have discussed when it can happen that $(t+1)d = 2n$.
Lemma~4.1 can also be used to construct other examples where the
quantity $td$ is fairly small when compared with $n$. For example,
suppose $G = PGL(2,p)$, where $p$ is a prime congruent to $3$
modulo $4$ and congruent to $\pm 1$ modulo $5$. Then $G$
has a maximal subgroup $H$ isomorphic to the symmetric group
$S_4$ and also a subgroup $K$ isomorphic to the alternating group
$A_5$ and such that $|H:H \cap K| = 2$. By Lemma~4.1, this yields
an example with $t = |K:H \cap K| = 5$ and $2n = 5r$, where
$n = |G:H|$ and $r = |G:K|$. By Lemma~4.2, we have $d \le r$, and
so $td \le 5r = 2n$.
\bigskip

\iitem{5. The case $t = n - 1$.}

In this section, we consider the case $t = n - 1$ to conclude our
analysis of equality in Theorem~C. Since we are assuming that
$(t + 1)d = 2n$, we see that $d = 2$ in this situation, and the
permutation module $F[S]$ has a $2$-dimensional submodule. The
following is an easy corollary of our Theorem~C.
\medskip

\iitem{(5.1) COROLLARY.}~~{\sl Let $S$ be a primitive right $G$-set
and suppose that $M \sbs F[S]$ is a $2$-dimensional submodule. If
$v \in M$ is nonzero, then either $t(v) = n$ or $t(v) = n - 1$, where
$n = |S|$. Also, if $t(v) = n - 1$, then $d(v) = 2$.}
\medskip

\iitem{Proof.}~~Write $t = t(v) > 0$ and note that $d = d(v) \le 2$
since $\gen v \sbs M$. By Theorem~C, we have
$2(t+1) \ge (t+1)d \ge 2n$, and thus $t + 1 \ge n$, as required. Also,
if $t = n - 1$, we have $d \ge 2$, and thus $d = 2$.\qed
\medskip

It is clear that in the situation of Corollary~5.1, some nonzero
element $v \in M$ satisfies $t(v) < n$, and thus $t(v) = n - 1$ and
$d(v) = 2$. Whenever $F[S]$ has a $2$-dimensional submodule,
therefore, we are in a situation where equality holds in
Theorem~C. We proceed to classify the primitive permutation
groups $G$ for which the permutation module $F[S]$ has a
submodule of dimension $2$. (For the sake of brevity, we will not
carry out this classification for all fields $F$, but only for
sufficiently large fields of some given characteristic.)
\medskip

\iitem{(5.2) THEOREM.}~~{\sl Let $S$ be a faithful primitive right
$G$-set and let $H$ be the stabilizer of a point in the action of $G$
on $S$. Suppose that the permutation module $F[S]$ has a
submodule of dimension $2$, where $F$ is a field. Then $H$ is
cyclic and $G$ has a normal elementary abelian $q$-subgroup $E$
for some prime $q$. Also, $EH = G$ and $H$ acts faithfully and
irreducibly on $E$. Finally, if the prime $q$ is not the characteristic
of $F$, then $|E| = q$ and $|H| \le 2$.}
\medskip

\iitem{Proof.}~~Let $M \sbs F[S]$ be a $2$-dimensional
$G$-submodule. As we have observed, $M$ must contain a vector
$v$ with $t(v) = n - 1$, where $n = |S|$, and we write
$S - \supp v = \{x\}$. It follows that the stabilizer in $G$ of $v$ is
contained in $G_x$, and in particular, $\ker M \sbs G_x$. But $G_x$
contains no nontrivial normal subgroup of $G$, and thus $M$
is a faithful $G$-module.

Suppose now that some element $g \in G$ fixes two distinct points
$y,z \in S$. We argue that $g = 1$, and thus the action of $G$ on
$S$ is either regular or Frobenius. If $w \in M$ is arbitrary, then
$y,z \not\in \supp{w - wg}$, and thus $t(w - wg) \le n - 2$. By
Corollary~5.1, therefore, $w - wg = 0$, and thus $g$ acts trivially
on $M$. But $M$ is faithful, and thus $g = 1$, as claimed.

The action of $G$ on $S$ is primitive, and so if it is regular, then
$|G|$ is prime, and there is nothing further to prove. We can thus
assume that the action of $G$ on $S$ is Frobenius, and hence there
is a regular normal subgroup $E$ by Frobenius' theorem. Also, $E$ is
nilpotent by Thompson's theorem, and it follows by primitivity that
$E$ is an elementary abelian $q$-group for some prime $q$. Also,
$G = EH$ and $H$ acts faithfully and irreducibly on $E$.

If $q$ is not the characteristic of $F$, then by standard properties
of Frobenius groups, it follows that no faithful $G$-module can
have dimension smaller than $|H|$. Thus $2 = \dim(M) \ge |H|$,
and so $|H| = 2$. Since $H$ acts irreducibly on $E$, it follows that
$|E| = q$, as wanted.

Finally, suppose that $q$ is the characteristic of $F$. Since the
normal $q$-subgroup $E$ is not contained in $\ker M$, it follows
that $M$ is not a simple $G$-module, and in particular, $M$ is not
simple as an $H$-module. But $M$ is completely irreducible as an
$H$-module by Maschke's theorem, since $q$ does not divide
$|H|$. (This is because in a Frobenius group, the orders of the
kernel and the complement are always coprime.) It follows that the
restriction $M_H$ is the direct sum of two $1$-dimensional
submodules, and since $M$ is faithful, we see that $H$ is abelian.
An abelian Frobenius complement, however, is necessarily cyclic,
and this completes the proof.\qed
\medskip

All of the possibilities allowed by Theorem~5.2 can actually occur.
To see this, we consider first the situation where $q$ (the prime
divisor of the regular normal subgroup $E$) is not the characteristic
of $F$. Then $G$ is either cyclic of order $q$ or dihedral of order
$2q$ and $|S| = q$. In this case, $F[S]$ is the direct sum of a trivial
module and a module $V$ of dimension $q - 1$. If $F$ is large
enough to contain a primitive $q\th$ root of unity, it is clear that a
simple submodule of $V$ has dimension $1$ or $2$, and in either
case, $F[S]$ has a $2$-dimensional submodule.

The more interesting case is where $q$ is the characteristic of $F$.
In this situation, we have $G = HE$, where $E$ is an elementary
abelian $q$-group, $H$ is cyclic and $H$ acts faithfully and
irreducibly on $E$. We choose the field $F$ so that $|F| = |E|$,
and we note that there is a subgroup $A \sbs F^\times$ with
$H \cong A$. In this situation, we can identify $G$ with the group
described in the following lemma.
\medskip

\iitem{(5.3) LEMMA.}~~{\sl Let $F$ be any field. If $a,b \in F$, with
$a \ne 0$, let $g_{a,b}:F \to F$ be the affine linear map defined by
$x \mapsto ax + b$. Now let $A \sbs F^\times$ be a subgroup. Then
the set $G = \{g_{a,b} \mid a \in A, b \in F\}$ is a transitive
permutation group on $F$. Also, $G$ is primitive if no nonzero
proper subgroup of $F^+$ is invariant under multiplication by $A$.}
\medskip

\iitem{Proof.}~~We compute that
$(x)g_{a,b}g_{c,d} = (ac)x + (bc+d)$, and it follows easily that $G$
is a group of permutations of $F$. Since $g_{a,b}$ carries $0 \in F$
to $b$, and $b \in F$ is arbitrary, we see that the action of $G$ on
$F$ is transitive. 

The stabilizer in $G$ of the point $0 \in F$ is the subgroup
$G_0 = \{g_{a,0} \mid a \in A\}$. If $G_0 \sbs H \sbs G$, where $H$ is
a subgroup, we let $K = \{t \in F \mid g_{1,t} \in H\}$. Then $K$ is a
subgroup of $F^+$ and it is easy to see that $H = G$ if $K = F$ and
$H = G_0$ if $K = 0$. Also, if $b \in K$ and $a \in A$, then $ab \in
K$, and therefore, if no nonzero proper subgroup of $F^+$ is
invariant under multiplication by $A$, it follows that $G_0$ is
maximal in $G$, and so $G$ is primitive on $F$.\qed
\medskip

In our situation, $H$ acts irreducibly on $E$, and so in the language
of Lemma~5.3, the action of $G$ on $F$ is primitive, and we take
$S = F$. We show in this situation that the permutation module
$F[S]$ has a $2$-dimensional $G$-submodule.

To avoid confusion, we distinguish $F$ from $S$ by writing
$s_x \in S$ to denote the element corresponding to $x \in F$. Take
$v = \sum_{x \in F} x s_x \in F[S]$ and note that $t = t(v) = n-1$,
where $n = |S| = |F|$. Also, let $w \in F[S]$ be the sum of the
elements of $S$. We compute that
$$
vg_{a,b} = \sum_{x \in F} x s_{ax + b} =
\sum_{y \in F} {y-b \over a} s_y = {1 \over a}(v - bw) \,,
$$
and thus $\gen v$ is contained in the subspace spanned by $v$ and
$w$. Also, the above calculation shows that $w \in \gen v$, and since
it is clear that $v$ and $w$ are linearly independent, we have
$d = 2$, and thus $\gen v$ is the desired $2$-dimensional submodule.
\bigskip

\iitem{6. Prime degree.}

In this section we establish Theorem~E, which asserts the
inequality $t + d > p$ in the case where $|S| = p$ is prime and the
field $F$ has characteristic $0$. As usual, $v \in F[S]$ is nonzero,
$t = t(v) = |\supp v|$ and $d = d(v) = \dimm F{\gen v}$.

We shall see that the inequality of Theorem~E is intimately
related to the following theorem of N.~G.~Chebotar\"ev.
\medskip

\iitem{(6.1) THEOREM (Chebotar\"ev).}~~{\sl Let $p$ be prime and
suppose that $\zeta \in \C$ is a primitive $p\th$ root of unity. Let
$V$ be the Vandermonde matrix with $(i,j)$-entry equal to
$\zeta^{ij}$, for $0 \le i,j \le p-1$. Then all square submatrices of
$V$ have nonzero determinant.}
\medskip

Chebotar\"ev's proof of this result is presented by P.~Stevenhagen
and H.~W.~Lenstra in their expository paper \ref\stelen. There are
also several other proofs of Chebotar\"ev's result in the literature,
and four of these are referenced in \ref\stelen. Yet another proof,
which we cannot resist mentioning, can be found in \ref\ie.

It is possible to interpret Theorem~6.1 in the spirit of the kinds of
inequalities we are considering in this paper. To see how to do this,
let $G = \gen z$ be a group of prime order $p$. The group of linear
characters $\hat G$ of $G$ is also cyclic of order $p$, and we
choose a generating character $\mu$. Write $\mu(z) = \zeta$ and
note that $\zeta$ is a primitive $p\th$ root of unity.

What would it mean to say that some square submatrix of the
Vandermonde matrix $[\zeta^{ij}]$ has determinant $0$? Clearly,
this is equivalent to the existence of subsets
$X,Y \sbs \{0,1,\ldots,p-1\}$, where $|X| = |Y|$, and
coefficients $a_x$ for $x \in X$, not all of them zero, such that
$\sum_{x \in X} a_x \zeta^{xy} = 0$ for all $y \in Y$. If we write
$v = \sum_{x \in X} a_x z^x \in \C[G]$, we see that this system of
equations can be rewritten as $\mu^y(v) = 0$ for $y \in Y$. In other
words, the failure of Chebotar\"ev's assertion would be equivalent
to the existence of some nonzero vector $v \in \C[G]$ such that
$\lambda(v) = 0$ for at least $t$ linear characters $\lambda$, where
$t = |\supp v|$.

We saw in Section~1 that if $v \in \C[A]$, where $A$ is an abelian
group, then the quantity $d = \dim(\gen v)$ is exactly equal to the
number of linear characters $\lambda$ of $A$ such that
$\lambda(v) \ne 0$, and so there are exactly $p - d$ linear
characters $\lambda$ of $A$ such that $\lambda(v) = 0$. The failure
of Chebotar\"ev's assertion, therefore, would be equivalent to the
existence of a nonzero vector $v \in \C[G]$ such that
$p - d \ge t$. In other words, Theorem~6.1 implies that the
inequality $t + d > p$ holds when $S = G$ has prime order $p$ and
$F = \C$. Conversely, we see that if we could find an independent
proof of this inequality, that would yield a proof of Chebotar\"ev's
result. (We will present just such a proof later in this section.)

Finally, we mention that the equivalence that we have just
established between the inequality $t + d > p$ and Chebotar\"ev's
assertion is valid for any field that contains a primitive $p\th$ root
of unity. As we shall see, however, our inequality can fail for some
such fields, and it follows that for those fields, the conclusion of
Chebotar\"ev's theorem is false.

To prove Theorem~E (assuming Chebotar\"ev's result) we reduce
the general problem to the case where $S = G$ has prime order $p$
and $F = \C$. We need the following easy observation.
\medskip

\iitem{(6.2) LEMMA.}~~{\sl Let $S$ be a finite set and suppose that
$F \sbs E$ are fields. Let $V \sbs F[S]$ be an $F$-subspace and write
$EV$ to denote the $E$-span of $V$ in $E[S]$. Then
$\dimm E{EV} = \dimm FV$.}
\medskip

\iitem{Proof.}~~An $F$-basis $\BB$ for $V$ clearly spans $EV$ over
$E$, and so it suffices to observe that $\BB$ is linearly independent
over $E$. This follows via elementary linear algebra, however,
since a homogeneous system of linear equations with coefficients
in $F$ that has a nontrivial solution over $E$ must also have a
nontrivial solution over $F$.\qed
\medskip

If $S$ is a $G$ set and $F \sbs E$ are fields, then a vector
$v \in F[S]$ can also be viewed as lying in $E[S]$, and it is obvious
that $t(v)$ does not change as we change our point of view from
the field $F$ to the field $E$. Since $d(v)$ is just the dimension of
the space spanned by the $G$-translates of $v$, we see by
Lemma~6.2 that $d(v)$ is also invariant under this change of field.
\medskip

\iitem{Proof of Theorem E.}~~The group $G$ acts transitively on
the set $S$, which has prime cardinality $p$, and as we have
remarked previously, it is no loss to assume that the action of $G$
on $S$ is faithful. Let $P$ be a Sylow $p$-subgroup of $G$ and note
that $|P| = p$.

The $G$-module $F[S]$ can also be viewed as a $P$-module, and we
observe that the $P$-submodule generated by $v$ is contained in
the $G$-submodule generated by $v$. If we replace $G$ by $P$,
therefore, the value of $d$ may decrease. Of course, this change
has no effect on $t$, and so to prove that $t + d > p$, it is no loss to
assume that $|G| = p$. In this case, we can assume that $S = G$ and
that the action of $G$ on $S$ is regular.

If we replace the given field $F$ by the subfield generated over the
rational numbers $\Q$ by the coefficients of $v$, this does not
change either $t$ or $d$. We can thus assume that $F$ is finitely
generated over $\Q$, and hence that $F \sbs \C$. We can thus
replace $F$ by $\C$ without changing $t$ or $d$. It therefore
suffices to prove the inequality $t + d > p$ in the case where
$G = S$ and $F = \C$. As we have seen, however, this case of
Theorem~E follows from Chebotar\"ev's theorem, and so the proof
is complete.\qed
\medskip

For the following discussion, fix a prime number $p$ and let $G = S$
have order $p$. By Theorem~E, we know that if the characteristic
of $F$ is $0$, then the inequality $t + d > p$ holds for all nonzero
vectors $v \in F[G]$. If $F$ has prime characteristic, however, this
inequality can fail. The following theorem gives a necessary and
sufficient condition for this failure to occur, where the condition is
expressed in terms of the polynomial ring $F[X]$. 

Let $z$ be a generator for $G$ and note that each vector
$v \in F[G]$ can be uniquely written in the form $f(z)$, where
$f \in F[X]$ and $\deg f < p$. The quantity $t = t(v)$ is exactly the
number of nonzero coefficients in the polynomial $f$, and we write
$t(f)$ to denote this number.
\medskip

\iitem{(6.3) THEOREM.}~~{\sl Let $G = \gen z$ be a group of prime
order $p$ and suppose that $v \in F[G]$ is nonzero, where $F$ is an
arbitrary field. Write $v = f(z)$, where $f \in F[X]$ and
$\deg f < p$. Then $t(v) + d(v) \le p$ if and only if $t(f) \le \deg h$,
where $h(X) = \gcd{X^p - 1}{f(X)}$.}
\medskip

Before we proceed with the proof of Theorem~6.3, we show how
this result can be used to find explicit examples where
$t + d \le p$. That will enable us to find examples where
Chebotar\"ev's theorem fails in prime characteristic.

First, factor the polynomial $X^p - 1$ in $F[X]$ in order to determine
its proper divisors. If we can find such a divisor $h(X)$ such that
$t(h) \le \deg h$, we are done: simply take $v = h(z)$. (Note that
the condition $t(h) \le \deg h$ says that the polynomial $h$ is
``missing a term". In other words, the coefficient of $X^i$ in $h(X)$
is $0$ for some exponent $i < \deg h$.) If no proper divisor of
$X^p - 1$ is missing a term, then consider multiples $f(X)$ of
divisors $h(X)$ of $X^p - 1$ such that $\deg f < p$. If it is possible to
find such a multiple with $t(f) \le \deg h$, then $v = f(z)$ will be
the desired example.

Using a computer algebra system, it is easy to find prime numbers
$p$ and finite fields $F$ such that the polynomial $X^p - 1 \in F[X]$
has a proper divisor that is missing a term. A few such examples are
given in the following table. Of course, once an example is found,
the field $F$ can be replaced by any larger field, and so we list only
``minimal" examples. 
$$
\eqalign{
p = 7 &:~~ F=GF(2)\cr
p = 11 &:~~ F=GF(3)\cr
p = 13 &:~~ F=GF(3),~GF(4),~GF(5)\cr
p = 17 &:~~F=GF(2),~GF(13)\cr
p = 19 &:~~F=GF(4),~GF(5),~GF(7)\cr}
$$

We mention one further example: $p = 11$ and $F = GF(5)$. In this
case, no proper divisor of $X^{11} - 1$ is missing a term. One such
divisor, however, is $h(X) = X^5 + 2X^4 + 4X^3 + X^2 + X + 4$, and if
we compute $f(X) = (X - 2)h(X)$, we find that
$f(X) = X^6 + 3X^3 + 4X^2 + 2X + 2$. Then $\deg f = 6 < p$ and
$t(f) = 5 \le \deg h$. As we have seen, this yields an example where
$t + d \le p$.
\medskip

\iitem{Proof of Theorem~6.3}~~Suppose that $\deg f < p$, where
$f(X) \in F[X]$, and let $h(X) = \gcd{X^p - 1}{f(X)}$. Assuming that
$t(f) \le \deg h$, we write $v = f(z)$ and we work to control $d(v)$. 

Write $X^p - 1 = h(X)k(X)$ and let $s = \deg k$, so that
$s = p - \deg h \ge p - \deg f > 0$. Then $X^p - 1$ divides
$f(X)k(X)$, and since $z^p = 1$ in $G$, we see that
$vk(z) = f(z)k(z) = 0$. Now let $M \sbs F[G]$ be the subspace
spanned by the set $\{vz^i \mid 0 \le i < s\}$ and note that
$\dim(M) \le s$ since the spanning set has cardinality $s$. Since $k$
has degree $s$ and $vk(z) = 0$, we see that $vz^s$ is a linear
combination of the vectors $vz^i$ with $0 \le i < s$, and thus
$vz^s \in M$. It follows that $Mz \sbs M$, and thus $M$ is a
$G$-submodule of $F[G]$, and in fact, $M = \gen v$. Then
$d(v) = \dim(M) \le s = p - \deg h$ and we have
$t(v) + d(v) \le t(f) + (p - \deg h) \le p$ since we are assuming that
$t(f) \le \deg h$.

Conversely, suppose that $0 \ne v \in F[G]$ and that
$t(v) + d(v) \le p$. Write $v = f(z)$, where $f(X) \in F[X]$ and
$\deg f < p$. Then $t(f) = t(v)$, and we need to show that
$t(f) \le \deg h$, where $h(X) = \gcd{X^p - 1}{f(X)}$.

Choose the integer $s > 0$ as large as possible such that the set
$\{vz^i \mid 0 \le i < s\}$ is linearly independent. Since this set is
contained in $\gen v$, we clearly have $s \le d(v)$. By the
maximality of $s$, the vector $vz^s$ is a linear combination of
vectors in our set, and so we can write $vk(z) = 0$ for some
polynomial $k[X] \in F[X]$ of degree $s$. Then $f(z)k(z) = 0$ and it
follows by the division algorithm that the polynomial $f(X)k(X)$ is a
multiple of $X^p - 1$. (This is  because no polynomial of smaller
degree can vanish at $z$.) Then $X^p - 1$ divides $h(X)k(K)$, and
thus
$$
t(v) + d(v) \le p \le \deg h + \deg k = \deg h + s \le \deg h + d(v) \,.
$$
We conclude that $t(f) = t(v) \le \deg h$, as desired.\qed
\medskip

As we have seen, if $F$ is a field containing a primitive $p\th$ root
of unity, then the conclusion of Chebotar\"ev's theorem over $F$ is
equivalent to the assertion that $t + d > p$ for all choices of
nonzero vectors $v \in F[S]$. It follows that examples where
$t + d \le p$, such as those we presented earlier, yield examples
where the conclusion of Chebotar\"ev's theorem fails. (We need to
take extensions of our minimal fields that are large enough to
contain a primitive $p\th$ root of unity, and this, of course, is
always possible if $p$ is different from the characteristic of $F$.)

For each prime $p$, there are only finitely many characteristics
where Chebotar\"ev can fail, and thus there are only finitely many
characteristics where examples such as those discussed above can
occur. To see why this is true, consider the determinants of all
square submatrices of the complex matrix $[\zeta^{ij}]$, as in
Theorem~6.1. These are algebraic integers, and they are nonzero
by Chebotar\"ev's theorem, and so their norms are nonzero rational
integers. It should be reasonably clear that the characteristics
where the conclusion of Chebotar\"ev's can fail are are exactly the
primes that divide at least one of these integers, and clearly, there
are just finitely many such primes.  

We have seen several examples of primes $p$ and fields $F$ for
which it is possible to find vectors $v \in F[S]$ such that
$t + d \le p$. But this cannot happen if the characteristic of $F$ is
the given prime number $p$. 
\medskip

\iitem{(6.4) THEOREM.}~~{\sl Let $S$ be a transitive $G$-set with
$|S| = p$, a prime number, and let $0 \ne v \in F[S]$, where $F$
has characteristic $p$. Then $t(v) + d(v) > p$.}
\medskip

\iitem{Proof.}~~As in the proof of Theorem~E, we can assume that
$G = \gen z$ has order $p$. By Theorem~6.3, it suffices to show
that for all nonzero polynomials $f(X) \in F[X]$ with $\deg f < p$,
we have $t(f) > \deg h$, where $h(X) = \gcd{X^p - 1}{f(X)}$. But $F$
has characteristic $p$, and thus $X^p - 1 = (X - 1)^p$, and we see
that $h(X) = (X - 1)^m$ for some integer $m$. Our goal, therefore,
is to prove that $t(f) > m$. This, however, is immediate from the
following general lemma.\qed
\medskip

\iitem{(6.5) LEMMA.}~~{\sl Let $0 \ne f(X) \in F[X]$, where $F$ is an
arbitrary field. In the case where $F$ has prime characteristic $p$,
assume in addition that $\deg f < p$. Suppose that $(X - 1)^m$
divides $f(X)$, where $m \ge 0$. Then $t(f) > m$.}
\medskip

\iitem{Proof.}~~The result is clearly true if $m = 0$, and so we
assume $m > 0$ and we proceed by induction on $m$. If the
polynomial $f(X)$ has zero constant term, let $g(X) = f(X)/X$.
Then $g(X)$ is also a polynomial divisible by $(X - 1)^m$ and it too
satisfies the degree upper bound in prime characteristic. Since
$t(f) = t(g)$, we can replace $f$ by $g$. If we apply this argument
repeatedly, we can assume that that $f$ has a nonzero constant
term.

Now $f$ is not a constant polynomial since it is divisible by $X - 1$,
and since in prime characteristic $p$, we have $\deg f < p$, it
follows that the formal derivative $k(X) = f'(X) \ne 0$.
Furthermore, since $f$ has nonzero constant term, we see that
$t(k) < t(f)$. Finally, we note that $k(X)$ is divisible by
$(X - 1)^{m-1}$, and so by the inductive hypothesis, $t(k) > m - 1$.
Then $t(f) \ge t(k) + 1 > m$, as required.\qed
\medskip

Next, we prove Theorem~E again, but this time, we avoid appealing
to Chebotar\"ev's theorem. By our remarks concerning the
relationship between these two results, this will yield a new proof
of Chebotar\"ev's theorem. We need the following lemma.
\medskip

\iitem{(6.6) LEMMA.}~~{\sl Let $S$ be a finite transitive right
$G$-set and let $0 \ne v \in \C[S]$. The following then hold.
\smallskip
\ritem{(a)} There exists an algebraic number field $K$ and a
nonzero vector $v_0 \in K[S]$ such that $t(v_0) = t(v)$ and
$d(v_0) \le d(v)$.
\smallskip
\ritem{(b)} Let $p$ be any prime. Then there exists a finite field
$F$ of characteristic $p$ and a nonzero vector $v_1 \in F[S]$ such
that $t(v_1) \le t(v)$ and $d(v_1) \le d(v)$.
\medskip}

To see the significance of this result, let $S$ be a transitive $G$-set
with $|S| = n$ and let $\Xi(x,y)$ be a real-valued function that is
monotonically increasing in each of its two real variables. Suppose
we want to prove that an inequality of the form $\Xi(t,d) > n$ holds
for all nonzero vectors $v \in \C[S]$, where $t = t(v)$ and
$d = d(v)$. By Lemma~6.6(b), it suffices to prove the same
inequality for all finite fields $F$ of characteristic $p$, where $p$
is some fixed prime.

We have already seen that to prove Theorem~E, we can assume
that $F = \C$, and so by Lemma~6.6(b), it is enough to establish
that the inequality $t + d > p$ always holds for finite fields of some
fixed characteristic $q$, and furthermore, we get to choose $q$.
But we have already done this; we can take $q = p$ by
Theorem~6.4. Theorem~E will thus follow once we prove
Lemma~6.6, and as we have seen, Chebotar\"ev's result will then
also follow. 

We remark that we really do need Theorem~6.4, despite the fact
that as we have seen, there are guaranteed to be infinitely many
characteristics for which the inequality $t + d > p$ holds. This is
because our proof that there were at most finitely many ``bad"
characteristics relied on Chebotar\"ev's theorem, which we are now
trying to prove.
\medskip

\iitem{Proof of Lemma~6.6.}~~Let $X$ be the $|G| \times |S|$
matrix over $\C$ in which the row corresponding to $g \in G$ is the
vector of length $|S|$ given by the coefficients of the translate
$vg$ of $v$, taken in some fixed order. Then $d = d(v)$ is the rank
of this matrix, and so if $e > d$, then every $e \times e$
submatrix of $X$ will have determinant $0$.

Now let $R$ be the $\Q$-subalgebra of $\C$ generated by the
$t = t(v)$ nonzero coefficients of $v$ and their reciprocals. Let $M$
be a maximal ideal of $R$ and note that the field $K = R/M$ is
finitely generated as an algebra over the image of $\Q$ in $K$,
and we identify this image with $\Q$. It follows by the Nullstellensatz
(see Theorem~30.8 of \ref\isa) that $K$ is a finite degree extension
of $\Q$. 

The images in $K$ of the nonzero coefficients of $v$ are invertible,
and hence they are nonzero. The image of $v$, therefore, is a
vector $v_0 \in K[S]$ such that $t(v_0) = t$. To compute $d(v_0)$,
we need to determine the rank of the appropriate $|G| \times |S|$
matrix $Y$ over $K$. But $X$ has entries in $R$ and $Y$ is the image
of $X$ in $K$. If if $e > d$, therefore, then all $e \times e$
submatrices of $Y$ have determinant $0$. The rank of $Y$ is thus at
most $d$, and we conclude that $d(v_0) \le d(v)$. This establishes
assertion (a).

Now $K$ is an algebraic number field, and we consider its ring $A$
of integers. Let $P \sbs A$ be a prime ideal containing $p$ and note
that $F = A/P$ is a finite field of characteristic $p$. We can replace
$v_0$ by a nonzero scalar multiple without affecting either
$d(v_0)$ or $t(v_0)$, and so by a standard fact from Dedekind
domain theory, we can assume that the coefficients of $v_0$ all lie
in $A$ but that not all of them lie in the ideal $P$. (See
Lemma~29.20 of \ref\isa.)

Let $v_1$ be the image of $v_0$ in $F[S]$. Then $v_1$ is nonzero
and $t(v_1) \le t(v_0)$. Also, $d(v_1)$ is the rank of a certain
$|G| \times |S|$ matrix $Z$ over $F$, and this matrix is the image
of a matrix over $A$ whose rank (over $K$) is $d(v_0)$. Reasoning
as before, we see that the determinant of every $e \times e$
submatrix of $Z$ is $0$ if $e > d(v_0)$. It follows that
$d(v_1) \le d(v_0)$ and the proof is complete.\qed
\bigskip

\iitem{7. Infinite groups.}

Are there any results for infinite groups that are analogous to our
inequalities? Suppose that $S$ is an infinite transitive right $G$-set
and $F$ is a field. We can, of course, continue to think about the
permutation module $F[S]$, which, by definition, consists entirely
of vectors with finite support. But $n$ is infinite (by assumption)
and if $v \in F[S]$ is nonzero, then $d$ is necessarily infinite, so
what can we hope to prove? Here is an easy result that suggests
that perhaps there is some theory here, although we have not
pursued it beyond this very special case.
\medskip

\iitem{(7.1) THEOREM.}~~{\sl Let $G = \gen z$ be an infinite cyclic
group acting regularly on the set $S$. Let $0 \ne v \in F[S]$, where
$F$ is an arbitrary field. Then the codimension $c$ of $\gen v$ in
$F[S]$ satisfies $t - 1 \le c < \infty$.}
\medskip

The inequality $c < \infty$ in Theorem~7.1 says that $\gen v$ is
``large", and this is consistent with the inequality $dt \ge n$ that
holds in the finite case. But in the finite case, $c = n - d$, and so
the inequality $t - 1 \le c$ in Theorem~7.1 would correspond to an
inequality of the form $t + d \le n + 1$, and this is in the direction
opposite  of what we might expect. (But compare this with
Lemma~F.) 
\medskip

\iitem{Proof of Theorem 7.1.}~~If $t(v) = 1$, then $\gen v = F[S]$,
and so $c = 0$ and we are done. We can assume, therefore, that
$t(v) > 1$. We prove that $c$ is finite by producing a finite subset
$X \sbs S$ such that $F[S] = F[X] + \gen v$, where, of course, the
subspace $F[X] \sbs F[S]$ is the linear span of $X$. For convenience,
we think of the set $S$ as a ``horizontal" linear array of points,
where the generator $z$ acts by a right shift of one unit.

Since $t(v) > 1$, it is possible to choose a nonempty subset
$X \sbs S$ consisting of consecutive points, and such that $X$ is
exactly one unit too short to contain any translate of $\supp v$.
Then the subspace $F[X] + \gen v$ contains the point just to the
right of the interval $X$ and also the point just to the left of this
interval. It follows that $F[X]z \sbs F[X] + \gen v$ and also
$F[X]z\inv \sbs F[X] + \gen v$. Then $F[X] + \gen v$ is mapped into
itself by $z$ and $z\inv$, and hence it is a $G$ submodule of
$F[S]$. We have seen, however, that $F[X] + \gen v$ contains a
point of $S$, and it follows that $F[X] + \gen v = F[S]$, as wanted,
and thus $c < \infty$.

Next, we argue that $F[X] \cap \gen v = 0$ by showing that if
$0 \ne w \in \gen v$, then $\supp w \not\sbs X$. We can write $w$
as sum $w = \sum_{i \in I} a_i vz^i$ for some finite subset $I \sbs \Z$,
where the coefficients $a_i \in F$ are nonzero. Let $r = \min(I)$
and $s = \max(I)$ and note the $\supp w$ contains the
leftmost point $x$ of $\supp vz^r$ and also the rightmost point $y$
of $\supp vz^s$. Since $r \le s$, however, the distance between
these points of $\supp w$ is at least the distance between the
leftmost and rightmost points of $\supp v$, and so by the choice of
$X$ it is not possible that $x$ and $y$ both lie in $X$. 

It follows now that $F[S] = F[X] \plusdot \gen v$, and so
$c = \dim(F[X]) = |X|$. By the choice of $X$, however, we see that
$t \le |X| + 1$, and so $t - 1 \le c$, as wanted.\qed
\medskip

We close this section by revisiting Rudio's lemma. We mentioned in
Section~2 that the strong form of Rudio's lemma does not hold in
general for infinite primitive groups. We now present an example
that demonstrates this.

Take $F = \R$, the real numbers, in Lemma~5.3, and let
$A \sbs {\R}^\times$ be the subgroup consisting of the positive real
numbers. It is clear that a subgroup of the additive group of $\R$
that is closed under multiplication by positive numbers will also be
closed under multiplication by all real numbers, and hence there
are no such subgroups other than $0$ and $\R$ itself. It follows that
the corresponding affine linear group $G$ acts primitively on $\R$.

Now let $X$ be the set of positive real numbers. It is easy to see
that the translate of $X$ under the group element
$x \mapsto ax + b$ is exactly the set $\{t \in \R \mid t \ge b\}$. (We
are using the assumption that $a > 0$ here, of course.) There is no
translate of $X$, therefore, that contains $0$ but not $1$, and so
the strong form of Rudio's lemma fails.
\bigskip

\iitem{8. Further remarks and questions.}

We certainly have not answered, or even attempted to answer, all
possible questions of the type considered in this paper. In this
section, we mention some other areas for possible future research
that have occurred to us.

Instead of limiting our attention to permutation modules, as we
have here, we could work in a somewhat more general context. We
shall say that a $G$ module $M$ over a field $F$ is {\bf monomial} if
it has a basis $S$ such that for each element $s \in S$ and each
group element $g \in G$, the translate $sg$ is a scalar multiple of a
member of $S$. (The permutation modules, therefore, are exactly
those monomial modules where all of these scalars are $1$.) In this
context, we replace the assumption that $S$ is a transitive $G$-set
with the requirement that the action of $G$ on the set
$\{Fs \mid s \in S\}$ is transitive. (We shall use somewhat
nonstandard language here, and refer to these as {\bf transitive}
monomial modules.)

We can still define $n = |S| = \dim(M)$ and also
$t = t(v) = |\supp v|$ and $d = d(v) = \dim(\gen v)$ for vectors
$v \in M$, and we can hope to prove results analogous to those we
have established for permutation modules. For example, the
inequality $td \ge n$ of Theorem~B continues to hold for transitive
monomial modules, and in fact, the proof of Theorem~B goes
through essentially unchanged. Perhaps some of our other results
would also generalize to monomial modules with little or no
change, but we will not pursue that here.

Another question concerns solvable primitive groups, and more
generally, primitive groups having a regular normal subgroup. We
saw that if equality holds in Theorem~C and $1 < t < n-1$, then the
group $G$ cannot have such a regular normal subgroup. In fact, in
almost all of the primitive examples we know, if $t$ and $d$ are
relatively small when compared with $n$, the group $G$ fails to
have a regular normal subgroup. (The exceptions here are when $t$
is very small and when it is very nearly equal to $n$.) This suggests
that substantially better inequalities might hold for primitive
groups that have regular normal subgroups, and in particular, for
solvable primitive groups. Perhaps in those cases, some additive
inequality should hold, such as that in Theorem~E, which asserts
that $t + d > n$ in the case where $n$ is prime.

Also, we know that equality can occur in Theorem~C in the cases
where $1 < t < n-1$, but only when $t = d$ and $F$ has
characteristic $2$. This suggests the possibility that substantially
better inequalities might hold when $|t - d|$ is large or for fields
of other characteristics.

In order to find good inequalities of the type we have been
discussing here, it would be useful to have more examples of
extreme or nearly extreme cases. Given a nonzero submodule $M$ of
the permutation module $F[S]$, we write
$t(M) = \min \{t(v) \mid 0 \ne v \in M\}$. By Lemma~F, we know that
$t(M) \le n + 1 - \dim(M)$, but what more can be said? In particular,
in the case where $F$ is a finite field, it would be pleasant if we
could study this question computationally, but we see no good
algorithm for finding nonzero vectors $v \in M$ with small support.

We can say a bit about the effect of a change of field on the
quantity $t(M)$. In fact, the following easy result has nothing to
do with modules, and so we extend the definition of $t(M)$ to
arbitrary vector subspaces $M$ of $F[S]$.
\medskip

\iitem{(8.1) LEMMA.}~~{\sl Let $S$ be a finite set and suppose that
$F \sbs E$ are fields. Let $V \sbs F[S]$ be a nonzero $F$-subspace
and write $EV$ to denote the $E$-span of $V$. Then
$t(EV) = t(V)$.}
\medskip

Note that this result is somewhat analogous to Lemma~6.2, which
was the nearly trivial observation that $\dim(EV) = \dim(V)$.
\medskip

\iitem{Proof of Lemma~8.1}~~Since $V \sbs EV$, the inequality
$t(EV) \le t(V)$ is obvious. To prove the reverse inequality, let $\BB$
be an $F$-basis for $V$ and let $w \in EV$ with $t(w) = t(EV)$. We
can write $w = \sum_{b \in \BB} e_b b$, where $e_b \in E$, and not
all of the coefficients $e_b$ are $0$. If $s \in S$, then the
coefficient of $s$ in $w$ is $\sum_{b \in B} e_b b_s$, where we have
written $b_s$ to denote the coefficient of $s$ in the basis vector
$b$. In particular, if $s \not\in \supp w$, we have $\sum e_b b_s = 0$.

If we think of the coefficients $e_b$ as unknowns, we see that we
have a homogeneous system of $|S| - t(w)$ linear equations in
these unknowns, and the coefficients $b_s$ of this system lie in $F$.
(Recall that $b_s \in F$ since $\BB \sbs V \sbs F[S]$.) This system has
a nontrivial solution over $E$, and it follows by elementary linear
algebra that there is also a nontrivial solution over $F$. We can
therefore find coefficients $c_b \in F$, not all $0$, and such that
$\sum_{b \in \BB} c_b b_s = 0$ for all points $s \in S - \supp w$. In
particular, if we write $v = \sum c_b b$, then $0 \ne v \in V$ and
$\supp v \sbs \supp w$. Thus $t(V) \le t(v) \le t(w) = t(EV)$, and the
result follows.\qed
\medskip

\iitem{(8.2) COROLLARY.}~~{\sl Let $S$ be a right $G$-set and
suppose that $F \sbs E$ are fields such that each $G$-submodule
$M$ of $E[S]$ is spanned over $E$ by $M \cap F[S]$. Then if
$v \in E[S]$ is nonzero, there exists a nonzero vector $v_0 \in F[S]$
such that $t(v_0) \le t(v)$ and $d(v_0) \le d(v)$.}
\medskip

\iitem{Proof.}~~Let $M = \gen v$ and write $N = M \cap F[S]$, so
that $M = EN$, by hypothesis. By Lemma~8.1, we can choose a
vector $w \in N$ such that $t(w) = t(N) = t(M) \le t(v)$, where the
inequality holds because $v \in M$. Also, since $w \in N$, we have
$d(w) \le \dimm FN = \dimm EM = d(v)$, and the result follows.\qed
\medskip

Finally, we discuss possible improvements of Lemma~6.6(a). Recall
that according to that result, if we are given a nonzero vector
$v \in \C[S]$, where $S$ is a right $G$-set, then it is possible to find
an algebraic number field $K$ and a nonzero vector $v_0 \in K[S]$
such that $t(v_0) \le t(v)$ and $d(v_0) \le d(v)$. Since we can go
from the complex numbers $\C$ down to an algebraic number field
without increasing either $t$ or $d$, it seems reasonable to ask just
how small a subfield we can take. If the permutation module
$\C[S]$ happens to be multiplicity free, we can use Corollary~8.2 to
get a reasonably nice answer.
\medskip

\iitem{(8.3) THEOREM.}~~{\sl Let $S$ be a right $G$-set and
suppose that the corresponding permutation character $\pi$ of $G$
is multiplicity free. Let $F \sbs \C$ be the subfield generated by the
values of all of the irreducible constituents of $\pi$. If $v \in \C[S]$
is nonzero, then there exists a nonzero vector $v_0 \in F[S]$ such
that $t(v_0) \le t(v)$ and $d(v_0) \le d(v)$.}
\medskip

\iitem{Proof.}~~Since $\pi$ is a permutation character, it can be
afforded by a rational representation, and thus for each irreducible
character $\chi \in \irr G$, the Schur index $m(\chi)$ over $\Q$
divides the multiplicity $[\pi,\chi]$. (See Corollary~10.2(c) of
\ref\isc.) By assumption, each of these multiplicities is at most $1$,
and thus $m(\chi) = 1$ for all irreducible constituents $\chi$ of
$\pi$. It follows that each such constituent is afforded by an
$F$-representation, and thus there exists an $FG$-module $M$
affording the character $\pi$ and such that $M$ is the direct sum of
submodules that afford distinct irreducible characters.

Since $M$ and $F[S]$ are $F[G]$-modules that afford the same
character, they must be isomorphic. (For example, this follows
from Problem~9.5 of \ref\isc.) It follows that we can write $F[G]$
as a direct sum of submodules $M_\chi$ affording the distinct
irreducible constituents $\chi$ of $\pi$. Since $M_\chi$ affords
$\chi \in \irr G$, it follows that the submodules $\C M_\chi$
of $\C[S]$ are distinct and pairwise nonisomorphic and that $\C[S]$
is the direct sum of these submodules. Every submodule of $\C[S]$,
therefore, is the sum of some of the $\C M_\chi$. In particular,
every submodule of $\C[S]$ is the $\C$-span of its intersection with
$F[S]$, and we are in the situation of Corollary~8.2. The proof is
now complete.\qed
\vfil\eject
\parindent = 0pt
\frenchspacing

\centerline{REFERENCES}
\medskip

{\bf 1.}~~W. Burnside, {\it Theory of groups of finite order}, 2nd
ed., Dover Publications, New York, 1955.
\smallskip

{\bf 2.}~~P. J. Cameron, On groups of degree $n$ and $n-1$, and
highly-symmetric edge colourings, J. London Math. Soc. (2) {\bf 9}
(1974/75) 385--391.
\smallskip

{\bf 3.}~~P. J. Cameron, Finite permutation groups and finite
simple groups, Bull. London Math. Soc., {\bf 13} (1981) 1--22.
\smallskip

{\bf 4.}~~D. L. Donoho and P. B. Stark, Uncertainty principles and
signal recovery, SIAM J. of Appl. Math. {\bf 49} (1989) 906--931.
\smallskip

{\bf 5.}~~R. J. Evans and I. M. Isaacs, Generalized Vandermonde
determinants and roots of unity of prime order, Proc. of Amer.
Math. Soc., {\bf 58} (1977) 51--54.
\smallskip

{\bf 6.}~~B. Huppert, {\it Endliche Gruppen I}, Springer, Berlin-New
York, 1967.
\smallskip

{\bf 7.}~~I. M. Isaacs, {\it Character theory of finite groups},
Dover, New York, 1994.
\smallskip

{\bf 8.}~~I. M. Isaacs, {\it Algebra: A graduate course},
Brooks/Cole, Pacific Grove, 1994.
\smallskip

{\bf 9.}~~F. Rudio, \"Uber primitive Gruppen,
Journal f\"ur reine u. angew. Math., {\bf 102} (1888) 1--8.
\smallskip

{\bf 10.}~~P. Stevenhagen and H. W. Lenstra Jr., Chebotar\"ev
and his density theorem, Math. Intell., {\bf 18} (1996) 26--37.
\smallskip

{\bf 11.}~~T. Tao, An uncertainty principle for cyclic groups of
prime order, WorldWideWeb preprint
http://xxx.arxiv.cornell.edu/pdf/math.CA/0308286
\smallskip

{\bf 12.}~~A. Terras, {\it Fourier analysis on finite groups and
applications} LMS Student Texts {\bf 43}, Cambridge Univ. Press,
Cambridge (1999).
\smallskip

{\bf 13.}~~H. Wielandt, {\it Permutation Groups}, Academic Press,
New York, 1964.

\bye